\documentclass[12pt,reqno]{article}

\usepackage[usenames]{color}
\usepackage{amssymb}
\usepackage{amsmath}
\usepackage{amsthm}
\usepackage{amsfonts}
\usepackage{amscd}
\usepackage{graphicx}

\usepackage[colorlinks=true,
linkcolor=webgreen,
filecolor=webbrown,
citecolor=webgreen]{hyperref}

\definecolor{webgreen}{rgb}{0,.5,0}
\definecolor{webbrown}{rgb}{.6,0,0}

\usepackage{color}
\usepackage{fullpage}
\usepackage{float}
\usepackage{maple}

\usepackage{graphics}
\usepackage{latexsym}
\usepackage{epsf}
\usepackage{breakurl}
\usepackage{subfigure}
\usepackage{caption}

\begin{document}

\numberwithin{equation}{section}

\theoremstyle{plain}
\newtheorem{theorem}{Theorem}
\newtheorem{cor}[theorem]{Corollary}
\newtheorem{lemma}[theorem]{Lemma}
\newtheorem{prop}[theorem]{Proposition}

\theoremstyle{definition}
\newtheorem{definition}[theorem]{Definition}
\newtheorem{example}[theorem]{Example}
\newtheorem{conjecture}[theorem]{Conjecture}

\theoremstyle{remark}
\newtheorem{rem}[theorem]{Remark}
\numberwithin{theorem}{subsection}
\today

\begin{center}
{On a Generalized Moment Integral containing Riemann's Zeta Function:\\ Analysis and Experiment}
\vskip 1cm
Michael Milgram,\\
Consulting Physicist, Geometrics Unlimited, Ltd.,\\Box 1484, Deep River, Ont., Canada, K0J1P0.\\
\href{mailto:email}{\tt  mike@geometrics-unlimited.com} \\
https://orcid.org/0000-0002-7987-0820\\
and\\
Roy J. Hughes,\\
School of Mathematics and Statistics,\\
University of New South Wales,\\
Kensington, 2033, Australia.\\
\href{mailto:email}{\tt roy.hughes@unsw.edu.au}\\
https://orcid.org/0009-0009-1834-0111\\

\end{center}
\vskip .2 in
\begin{abstract}
Here, we study both analytically and numerically, an integral $Z(\sigma,r)$ related to the mean value of a generalized moment of Riemann's zeta function. Analytically, we predict finite, but discontinuous  values and verify the prediction numerically, employing a modified form of Cesàro summation. Further, it is proven and verified numerically that for certain values of $\sigma$, the derivative function $Z^{\prime}(\sigma,n)$ equates to one generalized tine of the Dirac comb function without recourse to the use of limits, test functions or distributions. A surprising outcome of the numerical study arises from  the observation that the proper integral form of the derivative function is quasi-periodic, which in turn suggests a periodicity of the integrand. This possibility is also explored and it is found experimentally that zeta function values offset (shifted) over certain segments of the imaginary complex number line are moderately auto-correlated.  
\end{abstract}
\section{Introduction}
\subsection{Background}
In a previous work \cite{milgram2024extension}, relationships between different parametric instances of inverse Mellin transform integrals of the form

\begin{equation}
Z(\sigma,r)\equiv \int_{-\infty}^{\infty}\frac{\zeta(\sigma\pm it)\,r^{\sigma\pm it}}{(\sigma\pm it)}\,dt\,= 
-i \int_{\sigma -i\,\infty }^{\sigma +i\,\infty }{\zeta \left(v \right) r^{v-1}}d v,\hspace{20pt}\sigma>1,
\label{Zdef}
\end{equation} 
were studied by the use of an extended form of Glasser's Master Theorem \cite{GlasserMaster} for special values of the real variables $\sigma$ and $r$. In a second previous work \cite{Mdet}, similar integrals were studied by generalizing $r$ to become a complex variable $r\rightarrow r\,\exp(i\phi)$, where it was found, with $n\in\mathbb{N}$ and $r\in \Re$, that for certain values of $r=n$, the value of the derivative function 
\begin{equation}
Z^{\prime}(\sigma,r)\equiv \frac{\partial}{\partial\,r}Z(\sigma,r)
\label{ZdefD}
\end{equation}
was indeterminate, depending on how the point $r=n$ was approached as a function of $\phi$. In particular, it was found that if the point $r=n$ was approached from certain directions in the complex $r$ plane, the function $Z^{\prime}(\sigma,n)$ yielded a completely consistent family of finite integrals. Approached from a different direction, the function $Z^{\prime}(\sigma,n)$ diverged (i.e., was singular). Here we reconsider and resolve that issue by first considering $Z(\sigma,r)$ as a function of $r$ and show that it is discontinuous, that is
\begin{align}
Z(\sigma,r)_{|_{\sigma>1}}=\int_{-\infty}^{\infty}\frac{\zeta \! \left(\sigma+i\,t  \right) r^{\sigma+i\,t }}{\sigma+i\,t }\,d t&=2\pi\lfloor r\rfloor,\hspace{20pt} &\sigma>1, r\neq n,\\
&=\pi(2n-1),&\sigma>1, r=n.
\label{ZsigRdef}
\end{align}
By differentiating the above with respect to $r$, it is clear that
\begin{equation}
Z^{\prime}(\sigma,r\neq n)_{|_{\sigma>1}}\equiv\frac{\partial}{\partial\,r}Z(\sigma,r)=\frac{1}{r}\,\int_{-\infty}^{\infty}\zeta\! \left(\sigma+i\,t  \right) r^{\sigma+i\,t }d t
 = 0\,,\hspace{40pt} \sigma>1, r\neq n
\label{ZsigRDrv}
\end{equation}
and, by studying numerical approximations to the derivative function $Z^{\prime}(\sigma,r=n)_{|_{\sigma>1}}$ at the discontinuity, we find that
\begin{equation}
Z^{\prime}(\sigma,n)_{|_{\sigma>1}}= \int_{-\infty}^{\infty}\zeta\! \left(\sigma+i\,t  \right) n^{\sigma+i\,t }d t
 = \infty\,,\hspace{120pt} \sigma>1\,.
\label{ZprInf}
\end{equation}

Combining \eqref{ZsigRDrv} and \eqref{ZprInf} identifies, for any positive integer $n$,
\begin{equation}
Z^{\prime}(\sigma,r)\equiv\int_{-\infty}^{\infty}\zeta\! \left(\sigma+i\,t  \right) r^{\sigma+i\,t }d t
 = 2\pi\,r\,\delta(r-n)\,,\hspace{120pt} \sigma>1\,,
\label{Zdelta}
\end{equation}
as one tine of the Dirac comb function when considered only as a function of the variable $r$.\footnote{The Dirac comb function is the set of Dirac delta functions of unit separation.} 
\begin{rem}
For a proof of \eqref{Zdelta}, see Section \ref{sec:siggt1}.
\end{rem}

Since the main interest in such integrals focusses on the range $0<\sigma<1$, in this work we consider that region by applying analytic continuation from the region $\sigma>1$.  Because $\zeta(\sigma+it)$ varies in sign throughout its range and $|\zeta(\sigma+it)|/\log(t)$ is bounded when $\sigma>1,~t>2$, \cite[Corollary 1, page 184]{Edwards}, it is reasonable to expect integrals of the form \eqref{Zdef} to be numerically convergent when $\sigma>1$. However, representations with $\sigma<1$, which are not expected to always converge numerically, are traditionally defined and given meaning by analytic continuation. 

First, we study the special case $\sigma=3/2$ and then employ straightforward translation of contour integrals (or the Master Theorem equivalent) to obtain identities valid for $\sigma=1/2$. This will be found in Section \ref{sec:SpecCase}. In Section \ref{sec:GenCase}, we then examine cases with general values of $\sigma$ following the same methods. Acknowledging that we are dealing with functions on the very edge of tractability both in a numerical and analytic sense, the following Section \ref{sec:numtests} applies the Cesàro regularization technique (see Appendix \ref{sec:Primer}) to ascertain if the analytic results are consistent with numerical approximation. Surprisingly, they are in excellent accord, given that the functions are discontinuous (see Remark \ref{rem:false}). Therefore we are reasonably confident that our analytic results are valid - the numbers do not lie. An unexpected digression arises when the various Cesàro approximations are viewed graphically, since they suggest that proper integrals associated with \eqref{Zdelta} are periodic, as are the integrands themselves. This observation is studied in a further series of numerical experiments, showing in Section \ref{sec:Correlat}, that the integrand function $|\zeta(\sigma+it)|$ is moderately auto-correlated . Finally, our discoveries are summarized and discussed along with suggested generalizations and applications.

\subsection{Notation and Lemmas} \label{sec:NotLem}
\subsubsection{Notation}

Throughout the following, $0<r\in \Re$, $n,m\in \mathbb{N}$ always and $:=$ means symbolic replacement. Other symbols are real except if specified. We employ $\lfloor ... \rfloor$  and $\lceil ... \rceil$ to represent the {\it floor} and {\it ceiling} functions respectively and $\delta(x)$ is the Dirac delta function. The $m^{th}$ derivative of $\zeta(s)$ is written $\zeta^{(m)}(s)$. Both computer programs Maple \cite{Maple23} and Mathematica \cite{Math23} were used throughout and are individually cited where necessary.

\begin{rem} \label{rem:Mark1}
It is important to emphasize that, if $x\in\Re$, here we define {\it floor} and {\it ceiling} functions such that $$\lfloor x \rfloor \text { means the greatest integer less than {\bf but not equal to} } x, $$ and
$$\lceil x \rceil \text { means the smallest integer greater than {\bf but not equal to} } x, $$ 
as opposed to the common (e.g., Maple, Mathematica) usage $\lfloor n \rfloor=n$ and $\lceil n \rceil=n$. In other words, the {\it floor} and {\it ceiling} functions are open at their respective ends. This means that the limit endpoints of $\lfloor x \rfloor$ and $\lceil x \rceil$ are undefined and any identity containing these functions requires that the values of that identity must be carefully specified when $x=n$ if the identity is to be complete (e.g. \cite[Eqs. II.1(3), Eqs. II.1(4) and II.5(16)]{VanDP}. This is a consequence of, and flows naturally from, the forthcoming analysis, where we independently obtain the value of a function at a point of discontinuity, but not necessarily as the mean of the limit of its values as the discontinuity is approached from above and below.
\end{rem}

\subsubsection{Lemmas - specific to $\sigma=3/2$}

From \cite[Eqs. 3.723(2) and 3.723(4)] {G&R} with $ \,r\,\in\Re$ and $j\in\mathcal{N}$, we have
\begin{align} \nonumber
\int_{-\infty}^{\infty}\frac{\sin \! \left(t\,\ln \! \left(\frac{j}{ \,r\,}\right)\right) t}{t^{2}+\frac{9}{4}}d t
 &= \;{\pi}{\left(\frac{j}{ \,r\,}\right)^{-\frac{3}{2}}},\hspace{20pt} &\mathrm{if } ~j> \,r\,;\\ \nonumber
 &= \,-{\pi}{\left(\frac{j}{ \,r\,}\right)^{\frac{3}{2}}},\hspace{20pt}& \mathrm{if } ~j< \,r\,;\\
 &=\; 0\hspace{20pt}& \mathrm{if } ~j= \,r\,,r=n, 
\label{J2s}
\end{align}
and
\begin{align} \nonumber
\int_{-\infty}^{\infty}\frac{\cos \! \left(t\,\ln \! \left(\frac{j}{ \,r\,}\right)\right)}{t^{2}+\frac{9}{4}}d t
 &= \frac{2}{3} \,\pi\left(\frac{j}{ \,r\,}\right)^{-\frac{3}{2}}\hspace{20pt} \mathrm{if } ~j> \,r\,;\\ \nonumber
 &= \frac{2}{3}{\pi}{\left(\frac{j}{ \,r\,}\right)^{\frac{3}{2}}},\hspace{20pt} \mathrm{if } ~j< \,r\,;\\
 &=\frac{2}{3}{\pi}\hspace{20pt} \mathrm{if } ~j= \,r\,, \,r\,\in\mathcal{N}\,. 
\label{G2}
\end{align}
\subsubsection{Lemmas - the general form}
From the same source, more general forms of the same listed identities are
\begin{equation}
\int_{-\infty}^{\infty}\frac{t\,\sin \! \left(a\,t \right)}{\sigma^{2}+t^{2}}\,d t
 = \pi \,{\mathrm e}^{-a\,\sigma},\hspace{20pt}a>0,
\label{L1}
\end{equation}
and
\begin{equation}
\int_{-\infty}^{\infty}\frac{\cos \! \left(a\,t \right)}{\sigma^{2}+t^{2}}\,d t
 = \frac{\pi}{\sigma} \,{\mathrm e}^{-a\,\sigma},\hspace{20pt}a>0.
\label{L2}
\end{equation}

\section{The special cases $\sigma=3/2$ and $\sigma=1/2$} \label{sec:SpecCase}

\subsection{The Master Theorem}

Consider the function $F(t)$ defined, for $ \,r\,>0$, $ \,r\,\in\Re$, by
\begin{equation}
F(t)\equiv \frac{\zeta \! \left(\frac{1}{2}+i\,t \right)  \,r\,^{\frac{1}{2}+i\,t}}{\frac{1}{2}+i\,t}-\frac{\zeta \! \left(\frac{3}{2}-i\,t \right)  \,r\,^{\frac{3}{2}-i\,t}}{\frac{3}{2}-i\,t}.
\label{Fdef}
\end{equation}
It is easy to show that $F(t)+F(-i-t)=0$. Hence, from Glasser's Master theorem \cite{GlasserMaster}
\begin{equation}
\int_{-\infty}^{\infty}\left(\frac{\zeta \! \left(\frac{1}{2}+i\,t \right)  \,r\,^{\frac{1}{2}+i\,t}}{\frac{1}{2}+i\,t}-\frac{\zeta \! \left(\frac{3}{2}-i\,t \right)  \,r\,^{\frac{3}{2}-i\,t}}{\frac{3}{2}-i\,t}\right)d t
 = -2\,\pi \, \,r\,
\label{Ans}
\end{equation}
because only the residue at the point $t=-i/2$ contributes. 

\subsection{Evaluation: $\sigma=3/2$} \label{sec:method1}
We begin by considering the convergent integral
\begin{align} \nonumber
Z(3/2,r)&\equiv\,\int_{-\infty}^{\infty}\frac{\zeta \! \left(\frac{3}{2}-i\,t \right)  \,r\,^{\frac{3}{2}-i\,t}}{\frac{3}{2}-i\,t}d t\\
& = 
\moverset{\infty}{\munderset{j =1}{\sum}}\; {\left(\frac{j}{ \,r\,}\right)^{-\frac{3}{2}}}\int_{-\infty}^{\infty}{\mathrm e}^{i\,t\,\ln \left(\frac{j}{ \,r\,}\right)} \left(\frac{3}{2\,t^{2}+\frac{9}{2}}+\frac{i\,t}{t^{2}+\frac{9}{4}}\right)d t
\label{JJ1}
\end{align}
by writing
\begin{equation}
\zeta \! \left({3}/{2}-i\,t \right) = 
\moverset{\infty}{\munderset{j =1}{\sum}}\! \frac{1}{j^{\frac{3}{2}-i\,t}}
\label{Zid}
\end{equation}
and noting that the summation and integration can be transposed because both are convergent. When decomposed into its real and imaginary parts, we arrive at
\begin{align} \nonumber
\int_{-\infty}^{\infty}{\mathrm e}^{i\,t\,\ln \left({j}/{ \,r\,}\right)}&\left(\frac{3}{2\,t^{2}+\frac{9}{2}}+\frac{i\,t}{t^{2}+\frac{9}{4}}\right)d t \\
&= 
\int_{-\infty}^{\infty}\left(\left(\frac{3}{2}\right)\frac{\cos \! \left(t\,\ln \! \left(\frac{j}{ \,r\,}\right)\right)}{t^{2}+\frac{9}{4}}-\frac{t\,\sin \! \left(t\,\ln \! \left(\frac{j}{ \,r\,}\right)\right) }{t^{2}+\frac{9}{4}}\right)d t
\label{J}
\end{align}
noting that the imaginary parts of the integral vanish by antisymmetry over the integration range.

\subsubsection{Case: $ \,r\,=n$}

In the case that $ \,r\,$ is a positive integer $ \,r\,=n$, we employ \eqref{J2s} and \eqref{G2} to find
\begin{align} \nonumber
Z(3/2,n)&=\int_{-\infty}^{\infty}\frac{\zeta \! \left(\frac{3}{2}-i\,t \right) n^{\frac{3}{2}-i\,t}}{\frac{3}{2}-i\,t}d t
 \\  \nonumber
&=\pi\,\moverset{n -1}{\munderset{j =1}{\sum}}\! \left(\frac{j}{n}\right)^{-\frac{3}{2}} \left({\left(\frac{j}{n}\right)^{\frac{3}{2}}}+{\left(\frac{j}{n}\right)^{\frac{3}{2}}}\right)
+\moverset{n}{\munderset{j =n}{\sum}}\! \frac{n^{3}\,\pi}{j^{3}}+\pi\,\moverset{\infty}{\munderset{j =n +1}{\sum}}\! \frac{{\left(\frac{j}{n}\right)^{-\frac{3}{2}}}-{\left(\frac{j}{n}\right)^{-\frac{3}{2}}}}{\left(\frac{j}{n}\right)^{\frac{3}{2}}}\\ \nonumber
&=\pi\,\moverset{n -1}{\munderset{j =1}{\sum}}\; 2 + \pi +0 \\
&=\pi\,\left(2\,n-1\right).
\label{Q1}
\end{align}

\subsubsection{Case $ \,r\, \neq n $}
In the case that $ \,r\,>0$ is not a positive integer, that is $ \,r\,\neq n$, we again employ \eqref{J2s} and \eqref{G2} to find
\begin{align} \nonumber
Z(3/2,r) &=\int_{-\infty}^{\infty}\frac{\zeta \! \left(\frac{3}{2}-i\,t \right)  \,r\,^{\frac{3}{2}-i\,t}}{\frac{3}{2}-i\,t}d t\\ \nonumber
& = 
\pi\,\moverset{{\lfloor  \,r\, \rfloor}}{\munderset{j =1}{\sum}}\! \left(\frac{j}{ \,r\,}\right)^{-\frac{3}{2}} \left({\left(\frac{j}{ \,r\,}\right)^{\frac{3}{2}}}+{\left(\frac{j}{ \,r\,}\right)^{\frac{3}{2}}}\right)
+\pi\,\moverset{\infty}{\munderset{j ={\lceil  \,r\, \rceil}}{\sum}}\! \frac{{\left(\frac{j}{ \,r\,}\right)^{-\frac{3}{2}}}-{\left(\frac{j}{ \,r\,}\right)^{-\frac{3}{2}}}}{\left(\frac{j}{ \,r\,}\right)^{\frac{3}{2}}}\\ \nonumber
&=\pi\,\moverset{{\lfloor  \,r\, \rfloor}}{\munderset{j =1}{\sum}}2+0\\
&=2\,\pi\,\lfloor  \,r\,\rfloor\,.
\label{Q2}
\end{align}

\begin{figure}[h] 
\centering
\includegraphics[width=0.75\textwidth,height=.75\textwidth]{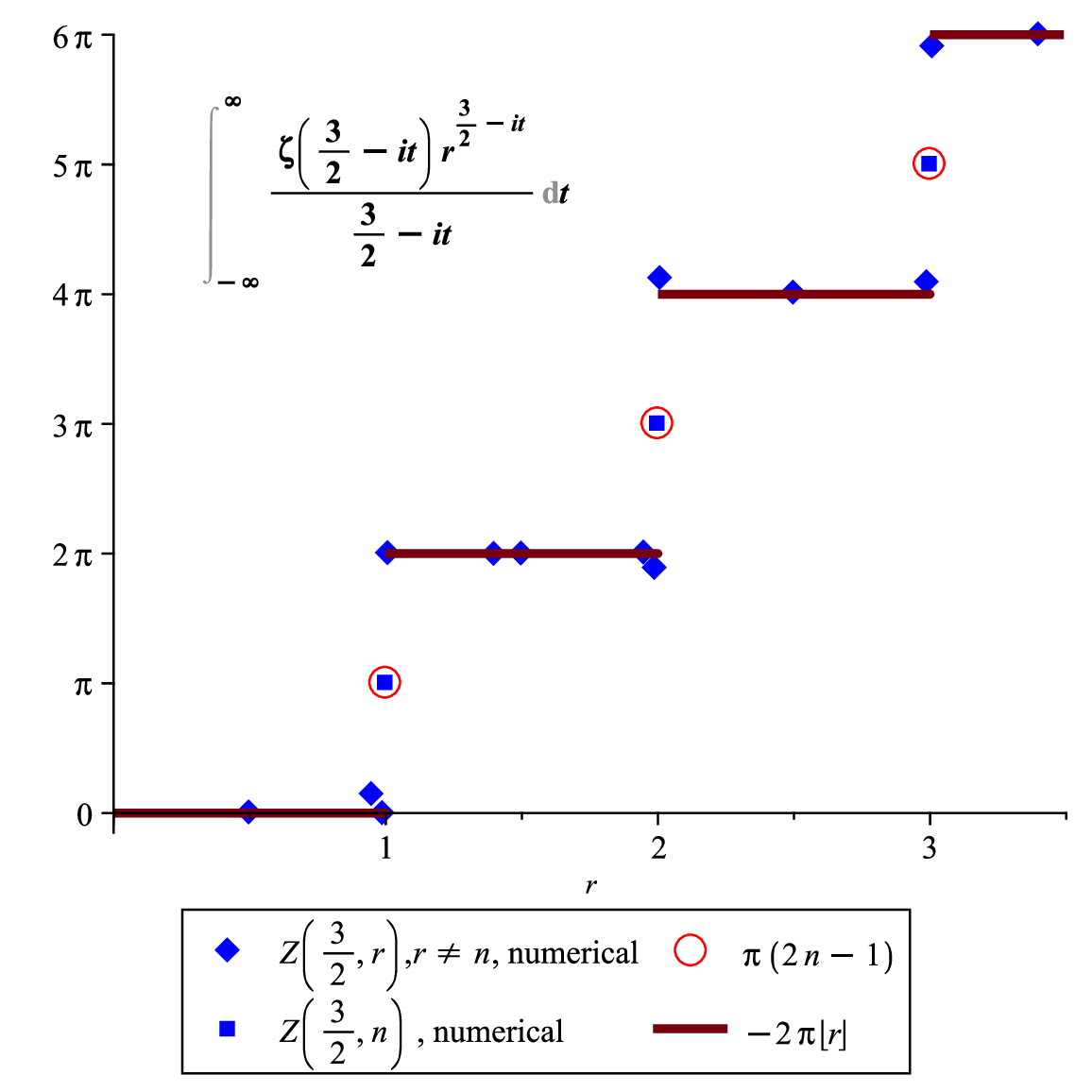}
\caption{The staircase function $Z(\frac{3}{2},r)$ over a small range of $r$, compared analytically and numerically. Note that the left and right limits of each of the horizontal ``treads" are open and that the values at the midpoints of the ``risers" are obtained by both an analytic and a numerical evaluation of the integral, not by decree.}
\label{fig:Z(r)}
\end{figure}

Figure \ref{fig:Z(r)} shows both the analytic results obtained above, as well as a few evaluations of Z(3/2,r) scattered over different values of $r$ obtained by (difficult) numerical integration, demonstrating substantial agreement.

\begin{rem} \label{sec:rgoeston}
Note that \eqref{Q2} does not reduce to \eqref{Q1} if $r\rightarrow n$ (see Remark \ref{rem:Mark1}) and therefore $Z(3/2,r)$ is discontinuous as a function of $r$ at $r=n$. However, at $r=n$, the function $Z(3/2,n)$ does lie exactly at the midpoint of the discontinuity, both analytically and numerically, as expected \cite[Eq. II.1(page 8)]{VanDP}.  
\end{rem}

\begin{rem} \label{rem:false}
As alluded to previously, the derivations presented here deal with subtle issues of limits and discontinuities. To emphasize this point, Appendix \ref{sec:false} presents homologous, but incorrect, identities obtained by a false derivation. The resolution is left as a challenge for the reader.
\end{rem}

\subsection{Analytic continuation to $\sigma=1/2$}

As noted in \cite{milgram2024extension}, by pairing integrals living inside the critical strip $0\leq\sigma\leq1$ with companions that are tractable and live outside, it becomes possible to evaluate the companion integral $Z(1/2,r)$ by applying either the Master Theorem or analytic continuation and compare with previous results. From \eqref{Ans},\eqref{Q1} and \eqref{Q2} we find
\begin{align} \label{Tint2a}
Z(1/2,r)\equiv\int_{-\infty}^{\infty}\frac{\zeta \! \left(\frac{1}{2}+i\,t \right)  \,r\,^{\frac{1}{2}+i\,t}}{\frac{1}{2}+i\,t}d t
 &=\, 2\,\pi  \left({ \lfloor  \,r\,\rfloor-r } \right)\,,\hspace{30pt}&r\neq n,\\
 &=-\pi\,,&r=n.
\label{Tint2}
\end{align}

so that, when $n=2$, we have,
\begin{equation}
\int_{-\infty}^{\infty}\frac{\zeta \! \left(\frac{1}{2}+i\,t \right) 2^{i\,t}}{\frac{1}{2}+i\,t}d t
 = -\frac{\pi}{\sqrt{2}}\,
\label{Tint}
\end{equation}
in agreement with \cite[Eq. (4.13)]{milgram2024extension}. Also, if $n=2$ we find
\begin{equation}
\int_{-\infty}^{\infty}\frac{\zeta \! \left(\frac{3}{2}-i\,t \right) 2^{-i\,t}}{\frac{3}{2}-i\,t}d t
 = \frac{3\,\pi \,\sqrt{2}}{4},
\label{Tint2A}
\end{equation}
in agreement with \cite[Eq.(4.12)]{milgram2024extension}.
\subsection{The derivative} \label{sec:Deriv}

\begin{figure}[h] 
\centering
\includegraphics[width=0.8\textwidth,height=.5\textwidth]{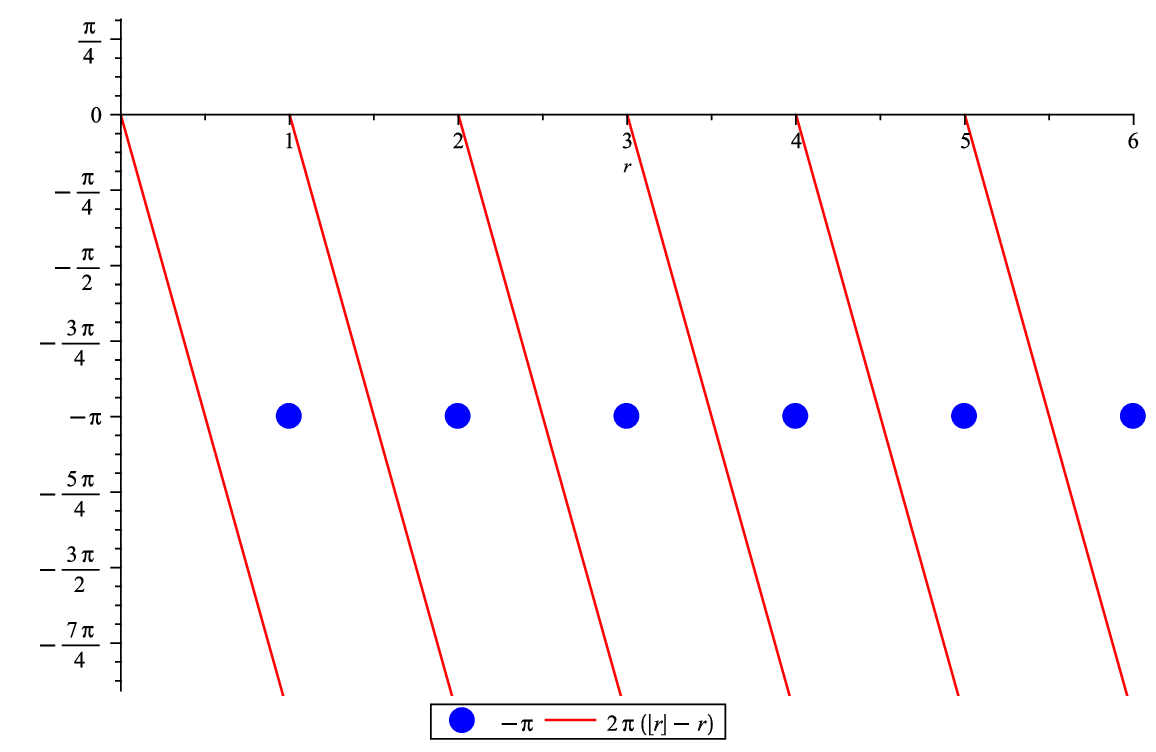}
\caption{The function $A(r)$ over a small range of $r$.}
\label{fig:Diagram}
\end{figure}

For typographical clarity, define the function 
\begin{align} \label{Rn}
A(r)&=2\,\pi(\lfloor r \rfloor -r),\hspace{30pt} ~ & r\neq n,\\
&=\,-\,\pi,&r=n,
\label{Rnn}
\end{align}
illustrated in Figure \ref{fig:Diagram}, demonstrating that $A(r)$ is invariant if $ \,r\,\rightarrow  \,r\,+m$ and that it can be characterized as a sawtooth function of period one and magnitude $2\pi$ where the special case $r=n$ corresponds to the midpoint of the sawtooth discontinuity. We also note that the slope of $A(r)$ is continuous as $ \,r\,\rightarrow n^{\pm}$, although the function itself is discontinuous, so it is relevant to query the value of the derivative $A^{\prime}(r)$ at $r=n$. By one simple, and standard definition,
\begin{align} \nonumber
A^{\prime}(n)\equiv\,A^{\prime}(r)\large|_{r=n}\equiv \underset{r\rightarrow n}{lim}\,\frac{d}{d\,r}A(r) &=\underset{r\rightarrow n}{lim}\,{\underset{h\rightarrow0}{lim}}\,\frac{A(r+h)-A(r)}{h}\\ \nonumber
&=2\pi\,\underset{r\rightarrow n}{lim}\,{\underset{h\rightarrow0}{lim}}\,\frac{\lfloor r+h \rfloor-r-h-\lfloor r \rfloor +r}{h} \\ \nonumber
&=2\pi\,\underset{r\rightarrow n}{lim}\,{\underset{h\rightarrow0}{lim}}\,\frac{\lfloor r \rfloor-h-\lfloor r \rfloor }{h}\\
&=-2\pi\,
\label{Ddef1}
\end{align}
and similarly if $h$ is replaced by $-h$, so the derivative $A^{\prime}(n)$ exists everywhere and is continuous according to this definition. However, by a second simple and standard definition

\begin{align} \nonumber
A^{\prime}(n)\equiv\,A^{\prime}(r)\large|_{r=n} &={\underset{h\rightarrow0}{lim}}\,\frac{A(n+h)-A(n)}{h}\\ \nonumber
&=2\pi\,{\underset{h\rightarrow0}{lim}}\,\frac{\lfloor n+h \rfloor-n-h +\pi/2}{h} \\ \nonumber
&=2\pi\,{\underset{h\rightarrow0}{lim}}\,\frac{\lfloor n \rfloor-n-h+\pi/2 }{h}\\
&=\infty
\label{Ddef2}
\end{align}
and the derivative is indeterminate when $r=n$, irrespective of our requirement that $\lfloor n \rfloor$ is undefined.

The important point here is that $A(r)$ represents the function $Z(1/2,r)$ (see \eqref{Tint2}) where the $r$ dependence only appears as an integrand term of the form $r^{1/2+it}$ on the left-hand side for which the derivative is well-known, continuous and consistent, independently of which of the two definitions are employed. Hence we have uncovered a pathology where the left-hand side of an identity always appears to be well-defined, and the right-hand side sometimes is not. Similar indeterminacies have been observed elsewhere \cite{Mdet} in related integrals, where it was speculated that there exists an associated, uncategorized, possibly essential, singularity at $r=n$, when the integral is studied as a function of complex $r$.
\label{sec:reduce}

\subsection{Differentiating}
Differentiating \eqref{Tint2} with respect to $r$ utilizing the definition \eqref{Ddef1}, yields
\begin{equation}
\int_{-\infty}^{\infty}\zeta \! \left({1}/{2}+i\,t \right)  \,r\,^{i\,t}d t
 = -2\,\pi \,\sqrt{ \,r}\,,\hspace{30pt} r\neq n,
\label{Q3}
\end{equation}
because
\begin{equation}
\frac{d}{d  \,r\,} \lfloor  \,r\, \rfloor =0\,,\hspace{100pt} r\neq n.
\end{equation}
Including the condition, \eqref{Q3} agrees with \cite[Eq. (5.28)]{Mdet}, which was obtained by an inverse Mellin transform, and reads
\begin{equation}
\int_{0}^{\infty}\Re\left(\zeta \! \left({1}/{2}+i\,t \right)  \,r\,^{i\,t}\right)d t
 = -\pi \,\sqrt{ \,r\,}\,,\hspace{30pt} r\neq n\,.
\label{Q3a}
\end{equation}
In the case that $r=n$, we also reproduce the indeterminacy alluded to above and first observed in \cite[Section (7.1)]{Mdet}. Utilizing \eqref{Ddef1} where $A^{\prime}(r)$ is continuous when $r\rightarrow n^{\pm}$ and no condition exists, by rewriting \eqref{Q3a} as
\begin{equation}
\int_{0}^{\infty}\left(\zeta_{I} \! \left({1}/{2}+i\,v \right)\,\sin \! \left(v\,\ln \! \left(r \right)\right)-\zeta_{R} \! \left({1}/{2}+i\,v \right)\,\cos \! \left(v\,\ln \! \left(r \right)\right)\right)d v
 = \pi \,\sqrt{r}
\label{Dr}
\end{equation} 
we reproduce \cite[Eq.(7.9)]{Mdet} when $r=n$; if we define the derivative $\frac{d}{dr}A(r)$ to be infinite at the point $r=n$ as in \eqref{Ddef2}, then we reproduce the indefinite value obtained in \cite[Eq.(7.10)]{Mdet} for that same integral and the right-hand side of \eqref{Dr} becomes indefinite. We reiterate that in this study we are considering a {\bf directed} limit $r\rightarrow n^{\pm}$ for real values of $r$, in contrast to \cite{Mdet} where the limit $r\rightarrow n$ was a {\bf directional} limit for complex values of $r$.  In \cite{Mdet}, it was shown that the finite choice represented by \eqref{Dr} (and accordingly \eqref{Ddef1}), leads to a consistent set of valuations for other integrals similar to \eqref{Dr}. Finally if $ \,r\,=1$, \eqref{Dr} demonstrates that
\begin{equation}
\int_{0}^{\infty}\zeta_{R} \! \left({1}/{2}+i\,t \right)d t
 = -\pi\,,
\label{Zr}
\end{equation}
in agreement with \cite[Eq. (5.26)]{Mdet}, although it is unlikely that \eqref{Zr} is numerically convergent and thereby it establishes one possible regularization of the integral among an infinite number of possibilites. 

\section{The General Case} \label{sec:GenCase}

\subsection{$\sigma>1$} \label{sec:siggt1}
As before, we begin by considering the convergent representation
\begin{equation}
\zeta \! \left(\sigma+i\,t  \right) = 
\moverset{\infty}{\munderset{j =1}{\sum}}\! \frac{1}{j^{\sigma}}\;{\mathrm e}^{-i\,t\,\ln \left(j \right)},\hspace{20pt}\sigma>1,
\label{Zsum}
\end{equation}
to be employed in the integrand of a convergent representation of $Z(\sigma,r)$, and then, in Section \ref{sec:SigLt1}, progress by analytic continuation to the region $\sigma<1$, where both the integral and/or the sum may not converge. Thus, with $\sigma>1$, 
\begin{equation}
Z(\sigma,r)\equiv\int_{-\infty}^{\infty}\frac{\zeta \! \left(\sigma+i\,t  \right) r^{\sigma+i\,t }}{\sigma+i\,t }\,d t=r^{\sigma} \moverset{\infty}{\munderset{j =1}{\sum}}\! \frac{1}{j^{\sigma}}\int_{-\infty}^{\infty}\left(\frac{\sigma\,\cos \left(t\,\ln \left(\frac{r}{j}\right)\right) }{\sigma^{2}+t^{2}}+\frac{\sin \left(t\,\ln \left(\frac{r}{j}\right)\right) t}{\sigma^{2}+t^{2}}\right)d t\,,
\label{Zgen}
\end{equation}
where the imaginary components vanish due to asymmetry. From \eqref{L1} it follows that
\begin{align}
\int_{-\infty}^{\infty}\frac{\sin \! \left(t\,\ln \! \left({r}/{j}\right)\right)}{\sigma^{2}+t^{2}}\,t\,d t
 &= \pi  \left(\frac{r}{j}\right)^{-\sigma},\hspace{20pt}&j<r; \label{J3a} \\
 &= \,-\pi  \left(\frac{r}{j}\right)^{\sigma},&j>r;
 \label{J4a}\\
 &=0, &j=r.
\end{align}
Similarly, from \eqref{L2} 
\begin{align}
\sigma  \int_{-\infty}^{\infty}\frac{\cos \! \left(t\,\ln \! \left(\frac{r}{j}\right)\right)}{\sigma^{2}+t^{2}}d t 
& = \pi  \left(\frac{r}{j}\right)^{-\sigma},\hspace{20pt}&j<r,
\label{J3b}\\
&=\pi  \left(\frac{r}{j}\right)^{\sigma},&j>r, \label{J4b}\\
&=\pi,& j=r. \label{J5b}
\end{align}
\subsubsection{$r\neq n$} \label{sec:RneqnSgt1}
If $r\neq n$, we find, as before (see \eqref{Q2})
\begin{equation}
Z(\sigma,r)=\int_{-\infty}^{\infty}\frac{\zeta \! \left(\sigma+i\,t  \right) r^{\sigma+i\,t }}{\sigma+i\,t }\,d t=2\pi\lfloor r\rfloor,\hspace{20pt} \sigma>1.
\label{Jt}
\end{equation}

Since $\sigma$ is now a continuous variable, we can consider differentiating \eqref{Jt} and \eqref{Jtn} with respect to both $r$ and $\sigma$. First, differentiating with respect to $\sigma$ (or integrating by parts), immediately yields the identity

\begin{equation}
\int_{-\infty}^{\infty}\frac{\zeta^{\left(1\right)}\! \left(\sigma+i\,t  \right) r^{\sigma+i\,t }}{\sigma+i\,t }d t
 = 
\int_{-\infty}^{\infty}\frac{\zeta \! \left(\sigma+i\,t  \right) r^{\sigma+i\,t }}{\left(\sigma+i\,t  \right)^{2}}d t-2\,\pi \, {\lfloor r \rfloor}\,\ln \! \left(r \right)
\label{Jd1}
\end{equation}
followed by differentiation with respect to $r$ which informs us that
\begin{equation}
\int_{-\infty}^{\infty}\zeta^{\left(1\right)}\! \left(\sigma+i\,t \right) r^{\sigma+i\,t }d t
 = 0
\label{JdaX}
\end{equation}
because $\frac{d{\lfloor r \rfloor}}{dr}=0$ if $r\neq n$. The generalization is obvious:
%
%
\begin{equation}
\int_{-\infty}^{\infty}\zeta^{\left(m\right)}\! \left(\sigma+i\,t  \right) r^{\sigma+i\,t }d t
 = 0
\label{JdaY}
\end{equation}
where $m\in \mathbb{N}$. Furthermore, differentiating \eqref{Jt} with respect to $r$ predicts that

\begin{equation}
\frac{1}{r}\,\int_{-\infty}^{\infty}\zeta\! \left(\sigma+i\,t  \right) r^{\sigma+i\,t }d t
 = 0\,, \hspace{60pt} \sigma>1,~r\neq n,
 \label{Jtd}
\end{equation}
raising a question about the potential interpretation of \eqref{Jtd} when $r=n$, where the definition of ``derivative" is ambiguous. See Section \ref{sec:numtests} below.

\subsubsection{$r=n$}

Exactly as in \eqref{Q1}, we have
\begin{equation}
Z(\sigma,n)=\int_{-\infty}^{\infty}\frac{\zeta \! \left(\sigma+i\,t  \right) n^{\sigma+i\,t }}{\sigma+i\,t }\,d t=\pi\left(2\,n-1\right),\hspace{20pt} \sigma>1,
\label{Jtn}
\end{equation}
and again, differentiating with respect to $\sigma$, gives

\begin{equation}
\int_{-\infty}^{\infty}\frac{\zeta^{\left(1\right)}\! \left(\sigma+i\,t \right) n^{\sigma+i\,t }}{\sigma+i\,t }d t
 = 
\int_{-\infty}^{\infty}\frac{\zeta \! \left(\sigma+i\,t \right) n^{\sigma+i\,t }}{\left(\sigma+i\,t  \right)^{2}}\,d t- \pi  \left(2\,n -1\right)\,\ln \! \left(n \right)\,.
\label{Jdn}
\end{equation}
Notice that the case $n=1$ applied to \eqref{Jdn} produces the identity

\begin{equation}
\int_{-\infty}^{\infty}\frac{\zeta^{\left(1\right)}\! \left(\sigma+i\,t \right) }{\sigma+i\,t }d t
 = 
\int_{-\infty}^{\infty}\frac{\zeta \! \left(\sigma+i\,t \right) }{\left(\sigma+i\,t  \right)^{2}}\,d t,\hspace{30pt}\sigma>1.
\label{Jdn1}
\end{equation}
If \eqref{Jdn1} is compared to \eqref{Jd1}, it suggests that  the (questionable) limit $r\rightarrow 1$ applied to \eqref{Jd1}, is valid.

\subsubsection{A Proof of \eqref{Zdelta}}
%
%
When $\sigma>1$, in Section \ref{sec:numtests} it will be shown that the integral \eqref{Jtd} appears to be (numerically) divergent only if $r=n$, and since it vanishes otherwise, we can justifiably write

\begin{equation}
Z^{\prime}(\sigma,r)\equiv\,\int_{-\infty}^{\infty}\zeta\! \left(\sigma+i\,t  \right) r^{\sigma+i\,t }d t
 = 2\pi\,r\,\delta(r-n)\,,\hspace{80pt} \sigma>1\,,n\neq 0.
\label{Zdelta2}
\end{equation}
The question that then arises is: ``Is $2\pi$ the appropriate normalization?" 
\begin{rem}
It is important to understand that \eqref{Zdelta2} differs from common invocations of the Dirac delta function (for examples, see \cite[Chapter V]{VanDP}) because it does not involve limits, test functions or distributions. The approach to the limit $r\rightarrow n$ is continuous, open-ended and always vanishes. {\bf In discontinuous fashion, only at the limit point $r=n$ does the integral diverge, and  this appears to be numerically true.} 
\end{rem}

\begin{theorem}
\begin{equation}
\int_{-\infty}^{\infty}\zeta\! \left(\sigma+i\,v  \right) r^{\sigma+i\,v }d v
 = 2\pi\,r\,\delta(r-n)\,,\hspace{120pt} \sigma>1\,,n\in\mathbb{N},
\label{Thm1}
\end{equation}

\begin{proof}
By a simple change of variables, we have, in more conventional notation,
\begin{equation}
\int_{-\infty}^{\infty}\zeta \left(\sigma+i\,v \right) r^{\sigma+i\,v}d v
 = 
-i \int_{\sigma -i\,\infty }^{\sigma +i\,\infty}\zeta \left(s \right) r^{s}d s 
\label{Xform}
\end{equation}
which we identify \cite[Eq. II.4.(13b)]{VanDP} as a 2-sided Fourier transform of the form
\begin{equation}
h \left(t \right) = 
\frac{1}{2\,\pi\,i} \int_{\sigma -i\,\infty }^{\sigma +i\,\infty }\frac{{\mathrm e}^{s\,t}\,f \left(s \right)}{s}ds 
\label{HT}
\end{equation}
where $t=\ln(r)$, $h(\ln(r))=2\pi\, r\,\delta(r-n)$, and $f(s)=2\pi\,s\,\zeta(s)$. The inverse transform of \eqref{HT} is \cite[Eq. II.4.(13b)]{VanDP}
\begin{equation}
f \left(s \right) = 
s \int_{-\infty}^{\infty}{\mathrm e}^{-s\,t}\,h \left(t\right) d t,\hspace{80pt} \Re(s)>1,
\label{FP}
\end{equation}
so first apply the change of variables $t=\ln(r)$ to  find \cite{Math23}
\begin{equation}
f \left(s \right) = 
s \int_{0}^{\infty}r^{-1-s}\,h \left(\ln \left(r \right)\right)d r, 
\label{X2}
\end{equation}
and then identify the various components of \eqref{HT} to obtain (with $n\neq 0$),
\begin{align}
\zeta \left(s \right)& = 
\int_{0}^{\infty}r^{-s}\,\delta \left(n-r  \right)d r\\
&=\sum_{n=1}^{\infty}1/n^s, \hspace{80pt}\Re(s)>1\,,
\label{X4}
\end{align}
the primary definition of $\zeta(s)$.
\end{proof}
\end{theorem}

\begin{rem}
The identity \eqref{JdaY} along with \eqref{Thm1} presents an interesting example of a non-constant function, all of whose derivatives vanish. See \cite{Stack240026} for a discussion about this point.
\end{rem}

\subsection{The case $\sigma<1$} \label{sec:SigLt1}

As suggested in \cite{milgram2024extension}, we now treat any of the above entities as a contour integral over a line $\Re(v)=\sigma$ in the complex $v$-plane where $v=\sigma+it$ by rewriting any of the above in the more conventional form of the contour integral representation (see \eqref{Zdef})

\begin{equation}
\int_{-\infty}^{\infty}\frac{\zeta \left(\sigma+i\,t  \right) r^{\sigma+i\,t }}{\sigma+i\,t }d t
 = 
-i \int_{\sigma -i\,\infty }^{\sigma +i\,\infty }{\zeta \left(v \right) r^{v-1}}d v ,
\label{Jeq}
\end{equation}
and translate the contour (equivalent to invoking the Master Theorem) such that $\sigma<1$ by accounting for the residues of the integrand so transited. First, we encounter a pole at $v=1$ with residue $r$; secondly a pole at $v=0$ with residue $-1/2$ and consider various cases.

\subsubsection{$r\neq n$}
From the above, we obtain
\begin{align}
\int_{-\infty}^{\infty}\frac{\zeta \! \left(\sigma+i\,t \right) r^{\sigma+i\,t}}{\sigma+i\,t}d t
 &= 2\,\pi  \left({\lfloor r \rfloor}-r \right),\hspace{20pt}&0< \sigma<1,
\label{Rneqn2}\\
&= 2\,\pi  \left({\lfloor r \rfloor}-r \right)+\pi, &\sigma<0,
\label{Rneqn3}
\end{align}
\eqref{Rneqn2} being an exact analogue of \eqref{Tint2a} -- there is no $\sigma$ dependence on the right-hand side other than the specification of the range of applicability, and it is invariant if $r:=r+m$. In exact analogy to subsection \ref{sec:RneqnSgt1}, differentiating first with respect to $\sigma$, then with respect to $r$ in both of \eqref{Rneqn2} and \eqref{Rneqn3} produces
\begin{equation}
\int_{-\infty}^{\infty}\zeta^{\left(1\right)}\! \left(\sigma+i\,t \right) r^{\sigma+i\,t}d t
 = 2\,\pi \,r\,\ln \! \left(r \right), \hspace{20pt}\sigma<1\,.
\label{Rneqn2d}
\end{equation}
Again, as in \eqref{Q3}, by differentiating \eqref{Rneqn2} or \eqref{Rneqn3} with respect to $r$ we obtain

\begin{align} \nonumber
\int_{-\infty}^{\infty}\zeta \left(\sigma+i\,t \right) r^{\sigma+i\,t}d t
 &= -2\,\pi r \hspace{20pt}&\sigma<1\,,r\neq n,\\
 &=-\pi r & \sigma=1\,,r\neq n,
\label{Rneqnd}
\end{align}
in agreement with \eqref{Q3a} by utilizing half the appropriate residue when $\sigma=1$.
\begin{rem} \label{rem:Mark2}
Since $\zeta(\sigma+it)\approx -i/t+O(t^{0})$ if $\sigma=1$, the imaginary part of the integrand contains the pole at $t=0$ and the imaginary integral that contains this divergent term, vanishes by anti-symmetry, leaving the finite result \eqref{Rneqnd}.
\end{rem}

\subsubsection{$r=n$}
For the case $r=n$, from \eqref{Jtn} plus the residues discussed above, we have

\begin{align} \label{Reqn2} \nonumber
\int_{-\infty}^{\infty}\frac{\zeta \! \left(\sigma+i\,t \right) n^{\sigma+i\,t}}{\sigma+i\,t}\,d t
 &=\pi(2n-1)-2n\pi\\
 &= -\pi\hspace{20pt} & 0<\sigma<1; \\
 &=\,~0 &\sigma<0,
\label{Reqn3}
\end{align}
consistent with \cite[Eqs. (4.3), (4.4) and (4.7)]{milgram2024extension}, all of which correspond to the case $n=1$. Differentiating with respect to $\sigma$ gives

\begin{align}
\int_{-\infty}^{\infty}\frac{\zeta^{\left(1\right)}\! \left(\sigma+i\,t \right) n^{\sigma+i\,t}}{\sigma+i\,t}d t -\int_{-\infty}^{\infty}\frac{\zeta \! \left(\sigma+i\,t \right) n^{\sigma+i\,t}}{\left(\sigma+i\,t \right)^{2}}d t 
 &= \pi\,\ln \! \left(n \right), \hspace{20pt} &0<\sigma<1,
\label{Reqn2d}\\
&=0&\sigma<0. \label{Reqn3d}
\end{align}
As in \eqref{Jdn1}, setting $n=1$ in  \eqref{Reqn2d} along with \eqref{Reqn3d} gives the generalization 
\begin{equation}
\int_{-\infty}^{\infty}\frac{\zeta^{\left(1\right)}\! \left(\sigma+i\,t \right) }{\sigma+i\,t }d t
 = 
\int_{-\infty}^{\infty}\frac{\zeta \! \left(\sigma+i\,t \right) }{\left(\sigma+i\,t  \right)^{2}}\,d t,\hspace{20pt} \forall\, \sigma.
\label{Jdn2}
\end{equation}
Again, following Section \ref{sec:SpecCase}, the value of \eqref{Rneqnd} when $r=n$ remains unclear, since the definition of derivative at a point is indeterminate. These issues are studied numerically in the following Section.

\section{Numerical tests}\label{sec:numtests}

Each of the integrals studied here is inherently oscillatory and therefore difficult to evaluate numerically, particularly since some of them may diverge. A simple way to explore the properties of such entities is to transform each into a (infinite) summation by the elementary act of subdividing the integration range into a large number of small, equal parts. By studying the numerical convergence of the sum employing Cesàro summation (known to provide a means of regularizing divergent -- or difficult -- sums), we obtain a means of verifying, or at least increasing confidence, in the identities that have been developed here (see Appendix \ref{sec:Primer}). 

\subsection{$\sigma>1$} \label{sec:sigGt1}


\begin{figure}[ht] \label{fig:sig4} 
\centering
\begin{subfigure}
[Cesàro approximation to the function defined in \eqref{sig4} with $\sigma=r=4$, suggesting divergence of order $T/2$. ]
{
\includegraphics[width=0.4\textwidth,height=.4\textwidth]{{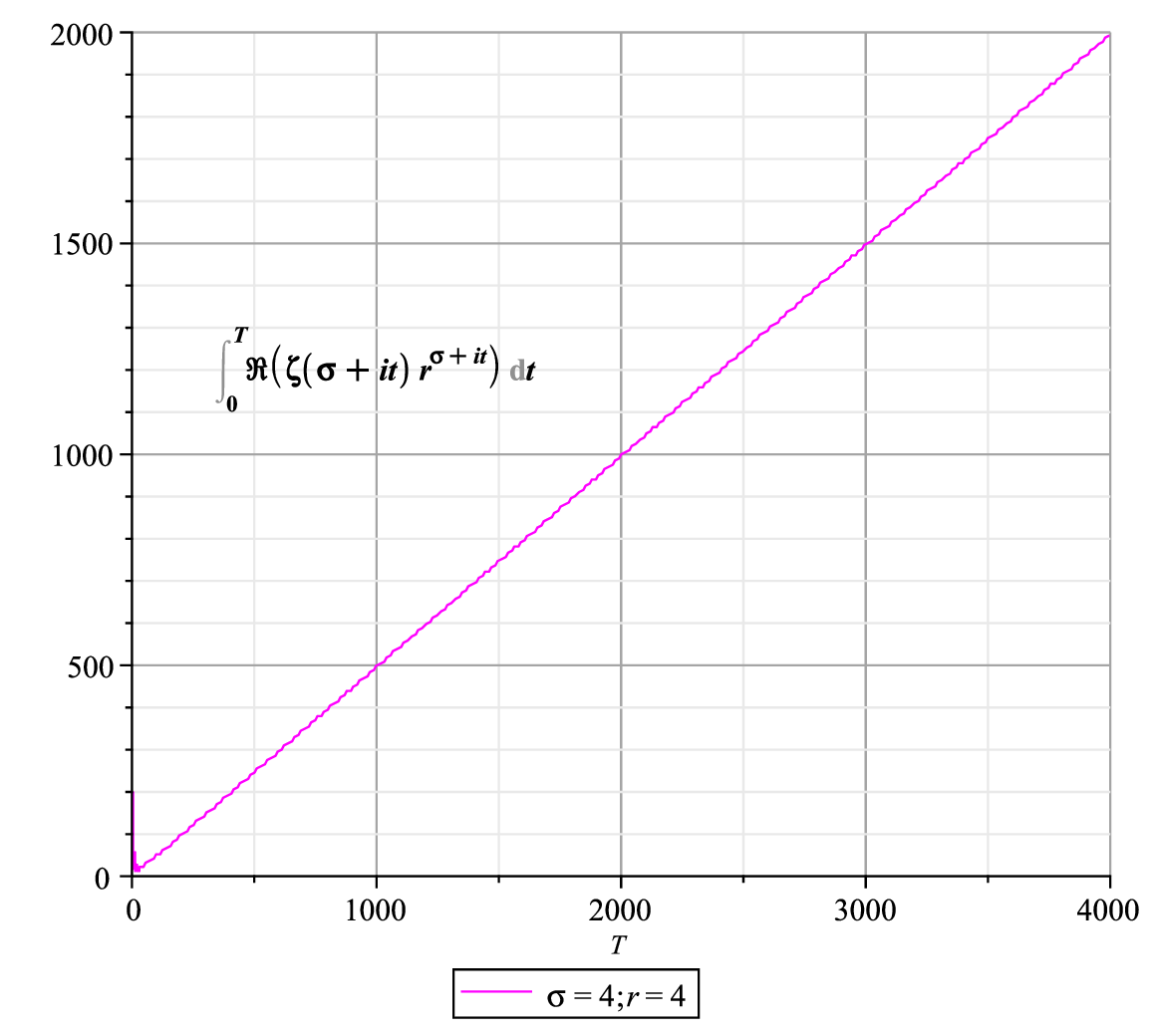}} \label{fig:Z1(r)}
}
\end{subfigure}
\hfill
\begin{subfigure}
[{Cesàro approximation to the function defined in \eqref{sig4} with $\sigma=4,~r=3.9$, suggesting convergence to zero.}] 
{
\includegraphics[width=0.4\textwidth,height=.4\textwidth]{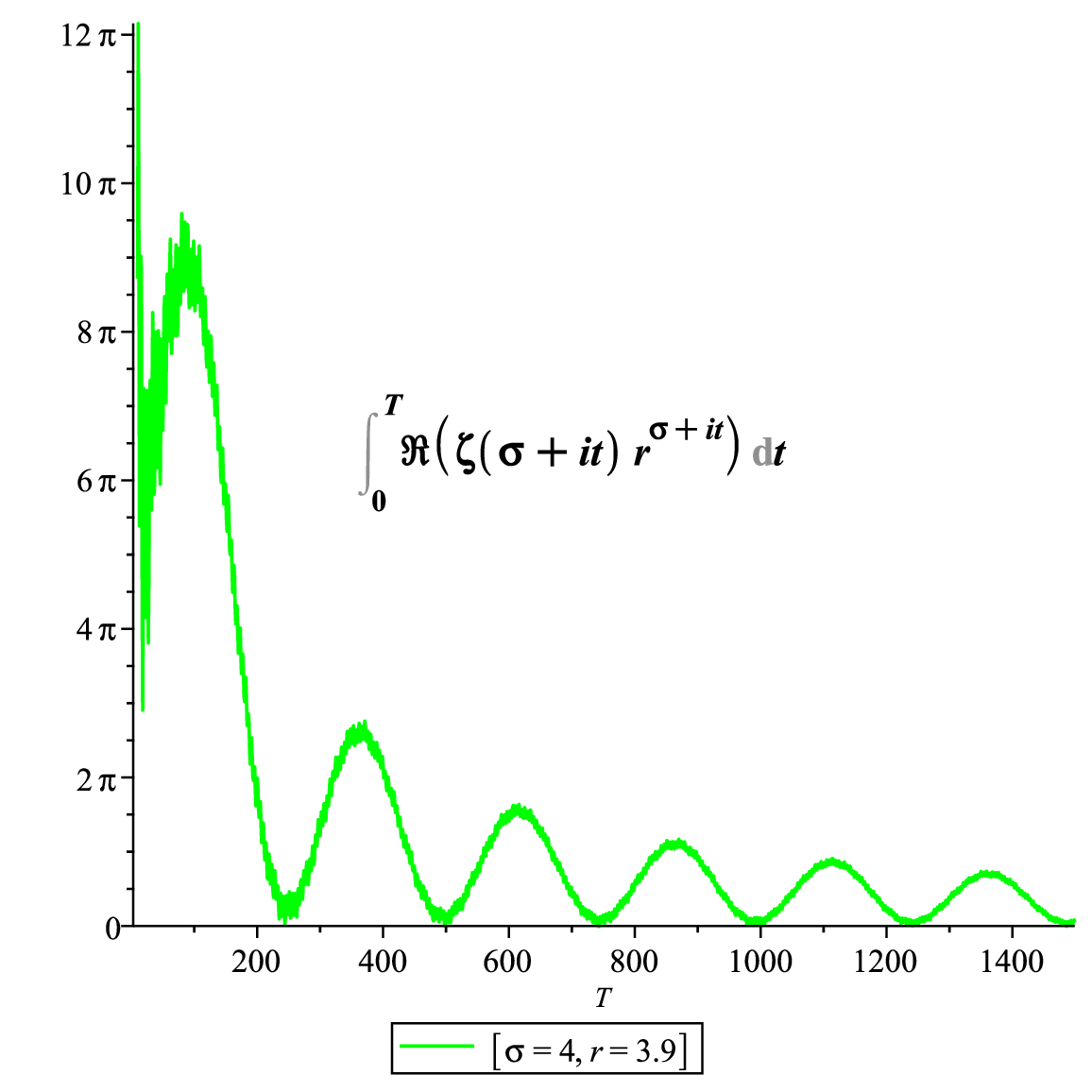} \label{fig:Z2(r)} 
}
\end{subfigure}
\caption{Numerical approximations to \eqref{sig4}}
\end{figure}

First, we consider the simple case $\sigma>>1$, (where the (oscillatory) integral is expected to converge when $r\neq n$ and there are no issues associated with analytic continuation), by setting $\sigma=4$ and comparing the cases $r=4$ and $r=3.9$. According to \eqref{Jtd}, the integral 
\begin{align}
{\underset{T\rightarrow\infty}{lim}}\int_{0}^{T}\Re\left(\zeta\! \left(4+i\,t  \right) r^{4+i\,t }\right)d t
 = 0\,, \hspace{20pt} r\neq n,
\label{sig4}
\end{align}
but could be either zero or infinite if $r=n$, depending on how the derivative is defined -- see Section \ref{sec:Deriv}. Figure (\ref{fig:Z1(r)}) is suggestive that \eqref{sig4} diverges to infinity of $\mathrm{O}(T^{1})$, at least when $r=4$. Similarly, when $r=3.9$, Figure (\ref{fig:Z2(r)}) is consistent with \eqref{sig4} and lends credibility to that identity, which states that the approximation should converge to zero.

%
%
\begin{figure}[h] \label{fig:Sighalf} 
\centering
\begin{subfigure}
[{Cesàro approximation to the function defined in \eqref{sighalf} with $\sigma=1/2,~r=0.9,~r=1.1$ and $r=2.1$, suggesting convergence to $-\pi r$ as T increases.}] 
{
\includegraphics[width=0.50\textwidth,height=.5\textwidth]{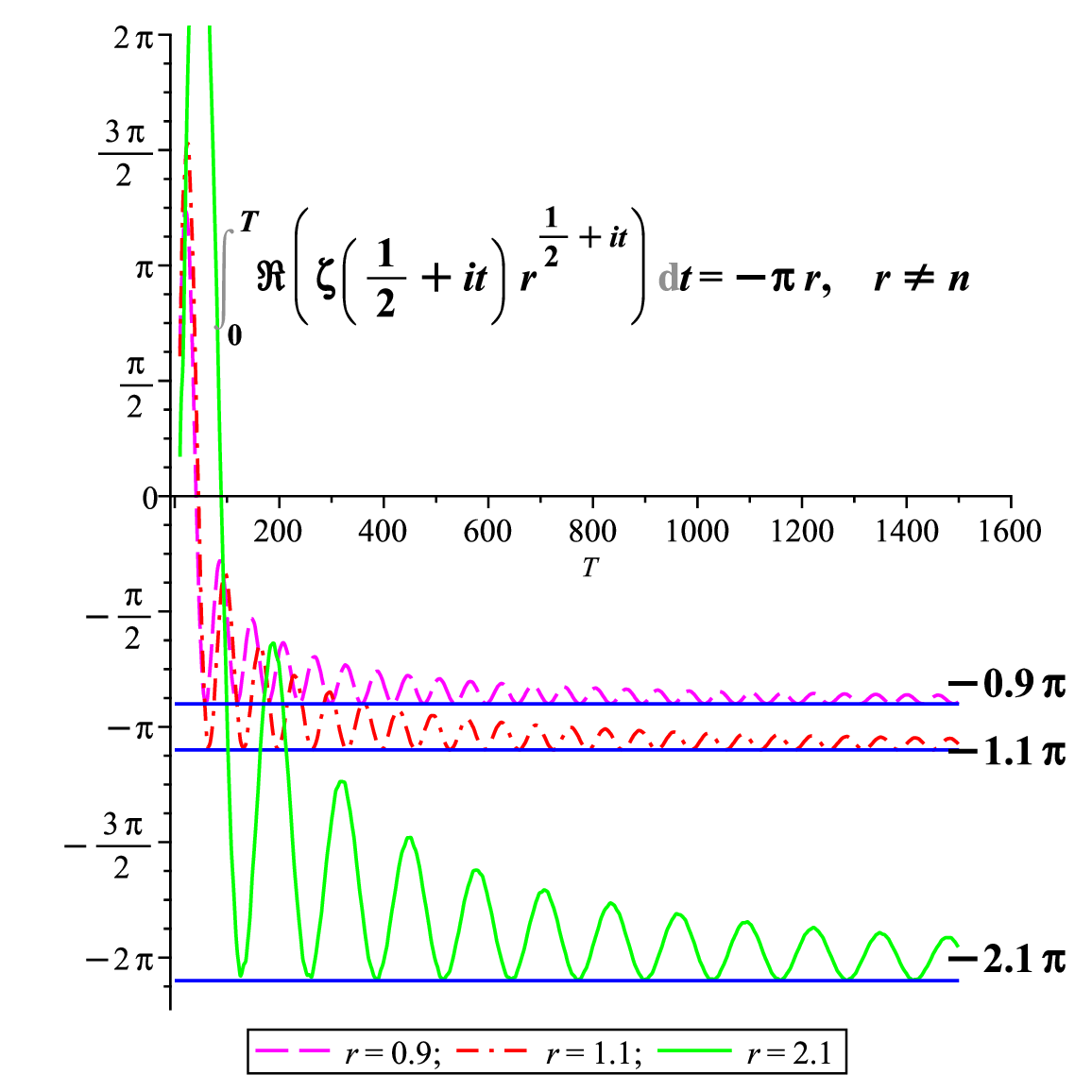} \label{fig:Z4(r)} 
}
\end{subfigure}
\hfill
\begin{subfigure}
[Cesàro approximation to the function defined in \eqref{sighalf} with $\sigma=1/2,~r=1$ and $r=2$, revealing a reasonable indication of divergence. The two cases are numerically indistinguishable within graphical resolution.]
{
\includegraphics[width=0.40\textwidth,height=.4\textwidth]{{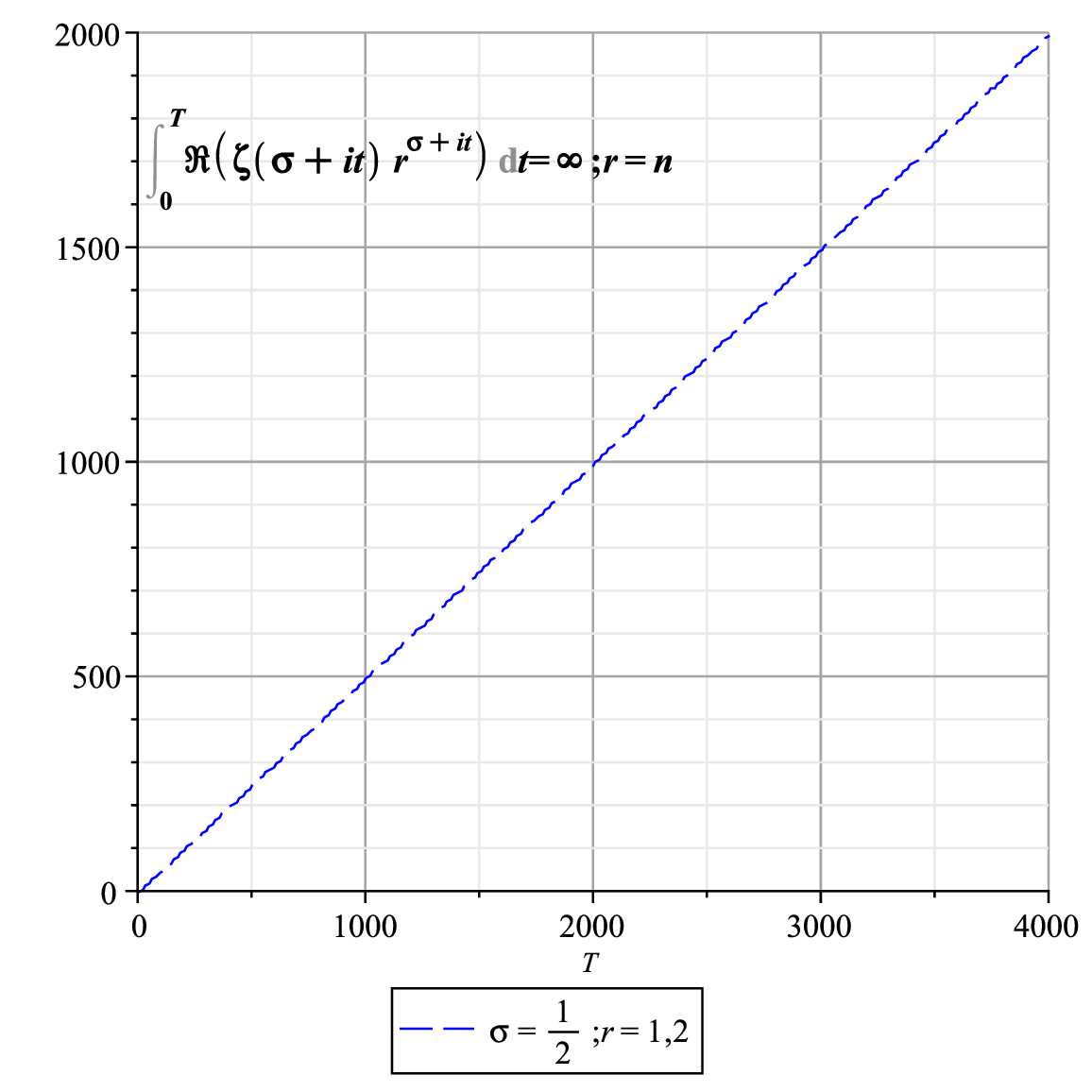}} \label{fig:Z3(r)}
}
\end{subfigure}
\caption{Numerical approximations to \eqref{sighalf}.}
\end{figure}

\subsection{$0<\sigma<1$} \label{sec:sigLt1}
It is also of interest to calculate the similar case where $\sigma<1$; that is, from \eqref{Rneqnd} we consider the identity

\begin{align}
\int_{0}^{\infty}\Re\left(\zeta\! \left(\sigma+i\,t  \right) r^{\sigma+i\,t }\right)d t
 = -\pi\,r \hspace{20pt} r\neq n,
\label{sighalf}
\end{align}
for several different values of $r$ and $\sigma$. Before doing so, we rewrite the finite version of \eqref{sighalf} in several equivalent forms:
\begin{align} \label{RFt}
\int_{0}^{T}\Re\left(\zeta \left(\sigma+i\,t  \right) r^{\sigma+i\,t}\right)& d t=\frac{r^{\sigma}}{2} \int_{-T}^{T}\zeta \left(\sigma-i\,t\, \right) {\mathrm e}^{-i\,t\,\ln \left(r \right)}d t  \\ 
&=r^{\sigma}\,\int_{0}^{T} \cos \left(t\,\ln \left(r \right)\right) \zeta_{R} \left(\sigma+i\,t \right)-\sin \left(t\,\ln \left(r \right)\right) \zeta_{I} \left(\sigma+i\,t \right)d t \label{RsinCos}\\
&=\int_{0}^{T}{| \zeta \left(\sigma +i\,t \right)|^{2}}\,\cos \left(t\,\ln \left(r \right)+\arg \left(\zeta \left(\sigma+i\,t  \right)\right)\right)\,dt,
\label{Rz}
\end{align}
and consider the above as $T\rightarrow \infty$. Inspired by Figure \ref{fig:Z2(r)}, we speculate that the integral itself, as a function of $T$, reflects the periodic nature of the trigonometric terms of the integrand. That this is more than speculation is illustrated in Figure \ref{fig:Z4(r)}, showing a slow, but steady convergence to $-\,\pi r$ using several values of $r\neq n$ and $\sigma=1/2$. In addition, the unexpected periodicity of the integral reappears, to be studied in Section \ref{sec:periodicity}. In the case $r=n$, Figure \ref{fig:Z3(r)} suggests that \eqref{sighalf} diverges of $\mathrm{O}(T^{1})$.

%
%
\begin{figure}[h]  
\centering
{
\includegraphics[width=1\textwidth,height=0.5\textwidth]{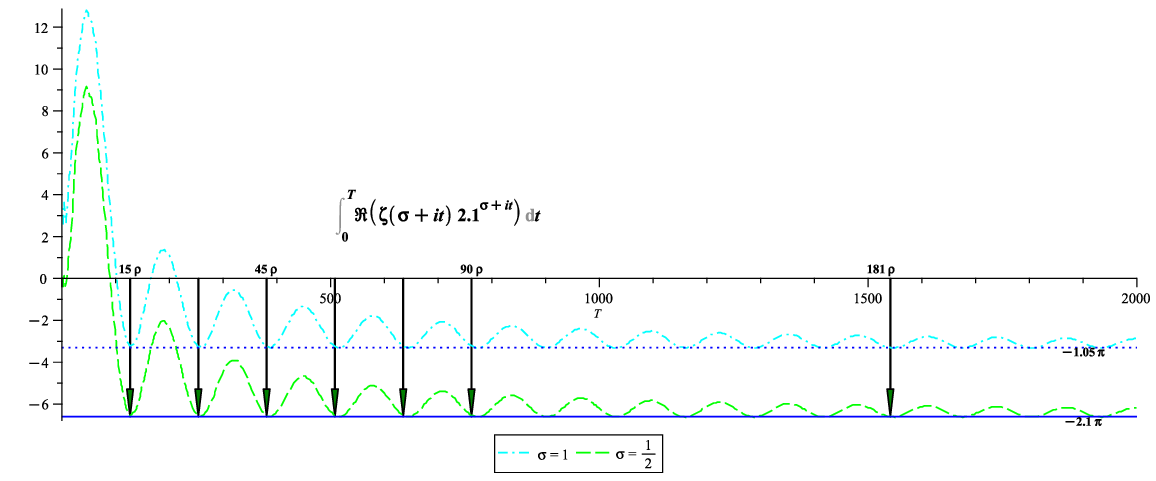}
} 
\caption{Cesàro approximations to the function defined in \eqref{sighalf} with $r=2.1$ and various $0<\sigma\leq 1$, with $\sigma=1/3,~ \sigma=3/4$ and $\sigma=7/8$ not shown because they are numerically indistinguishable, within the resolution of the figure, from the case $\sigma=1/2$ selected. The arrows mark integral multiples of a conjectured period $\rho$.}
\label{fig:3Sigs}
\end{figure}

When a similar calculation is performed for different values of $0<\sigma\leq 1$ and constant $r$, the Cesàro estimate also converges to the expected limit \eqref{sighalf}, which limit varies only if $\sigma=1$ as expected (see Remark \ref{rem:Mark2} and \eqref{Rneqnd}). This is illustrated in Figure \ref{fig:3Sigs}, which focusses on the Cesàro estimates for several values of $\sigma$ and constant $r=2.1$. This figure shows only the cases $\sigma=1/2$ and $\sigma=1$ because the others are indistinguishable from the case $\sigma=1/2$ within graphical resolution and again reveals a strong periodicity..

%
\begin{figure} [ht]
\noindent\begin{minipage} {0.47\textwidth} 
\centering
\includegraphics[width=1.0\textwidth,height=1.0\textwidth]{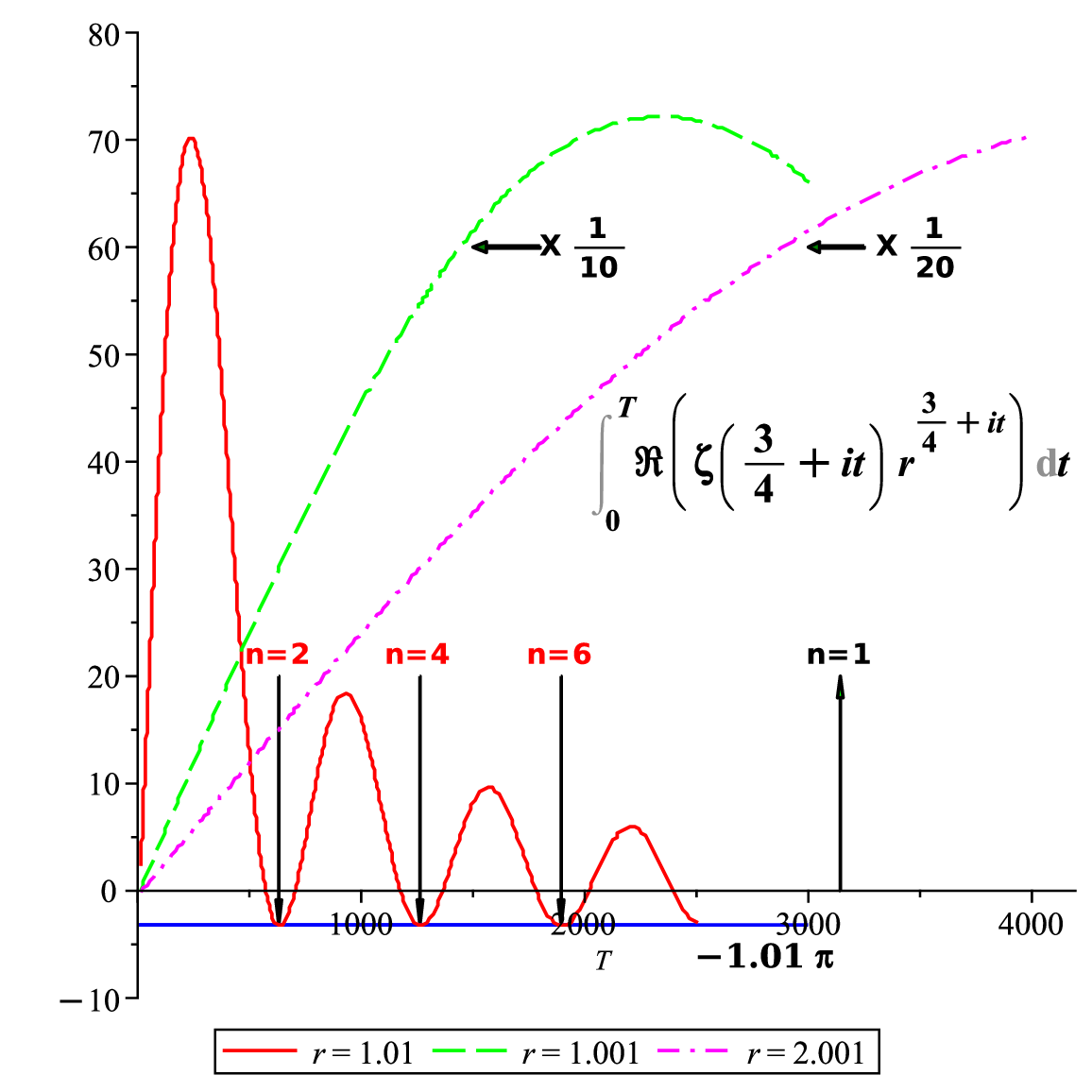} 
\captionof{figure}{Study with $r$ close to $n$. Down pointing vertical arrows indicate multiples of the solution points $t$ where $t \ln(1.01)=n\pi$, $n\geq 1$. Up pointing arrow labels two cases when $r\approx n$.}
\label{fig:twoRsNear1}
\end{minipage}%
\hfill
\begin{minipage} {0.47\textwidth} 
\centering
%
%
%
%

\includegraphics[width=1.0\textwidth,height=1.0\textwidth]{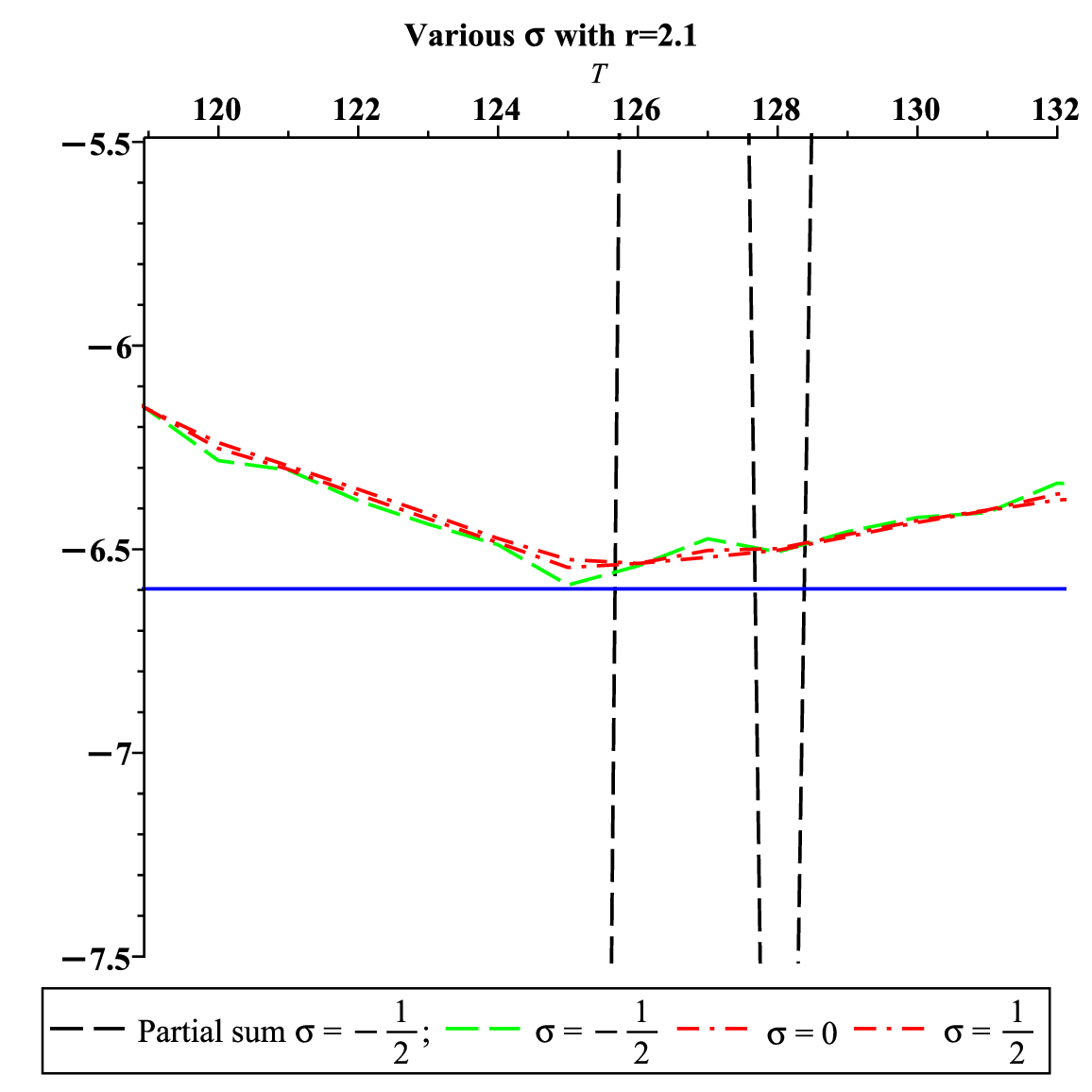}
\captionof{figure}{This is a detailed comparison between the Cesàro estimates with three values of $\sigma$ at constant $r=2.1$. It also shows that the partial sum for $\sigma=-1/2$ crosses the asymptotic line in three places.} 
\label{fig:MoreSigs2p1}
\end{minipage}%
\end{figure}

It is also of interest to question how \eqref{sighalf} converges to its discontinuous limit. Figure \ref{fig:twoRsNear1} shows that the convergence for $r=1.01$ and $r=1.001$ exists but it takes longer to ``turn over" as $r$ approaches unity. Significantly, once it does ``turn over", it will approach the finite limit $-\pi r$, rather than diverging. Since we shall shortly see that the oscillations in the partial sum appear to be influenced by a periodicity $\rho=2\pi/\ln(r)$, when $r=1$ the curve will never turn over, and \eqref{sighalf} becomes infinite at $r=1$. However, as Figure \ref{fig:twoRsNear1} shows, this does not lead to any insight concerning the case $r\approx 2$, or indeed any other integer -- see Figures \ref{fig:Z1(r)} and \ref{fig:Z3(r)}.

In a final study of this case, Figure \ref{fig:MoreSigs2p1} presents a detailed view near a minimum, of the Cesàro estimates with $\sigma=1/2$, $\sigma=0$ and $\sigma=-1/2$, all of which, according to \eqref{sighalf}, should converge to the same limit; as observed in this Figure they are effectively indistinguishable within graphical resolution. We also note that the partial sum  estimate for $\sigma=-1/2$ is far more varied than is the same estimate for the integral corresponding to $\sigma=1/2$ (shown below in Figure \ref{fig:FirstCross2p1} of Section \ref{sec:periodicity}), but it still intersects the asymptotic line $-2.1\pi$ near the minimum of the Cesàro estimates.

\begin{figure} [ht]
\noindent\begin{minipage} {0.40\textwidth} 
\centering
%
%
\includegraphics[width=1.0\textwidth,height=1.0\textwidth]{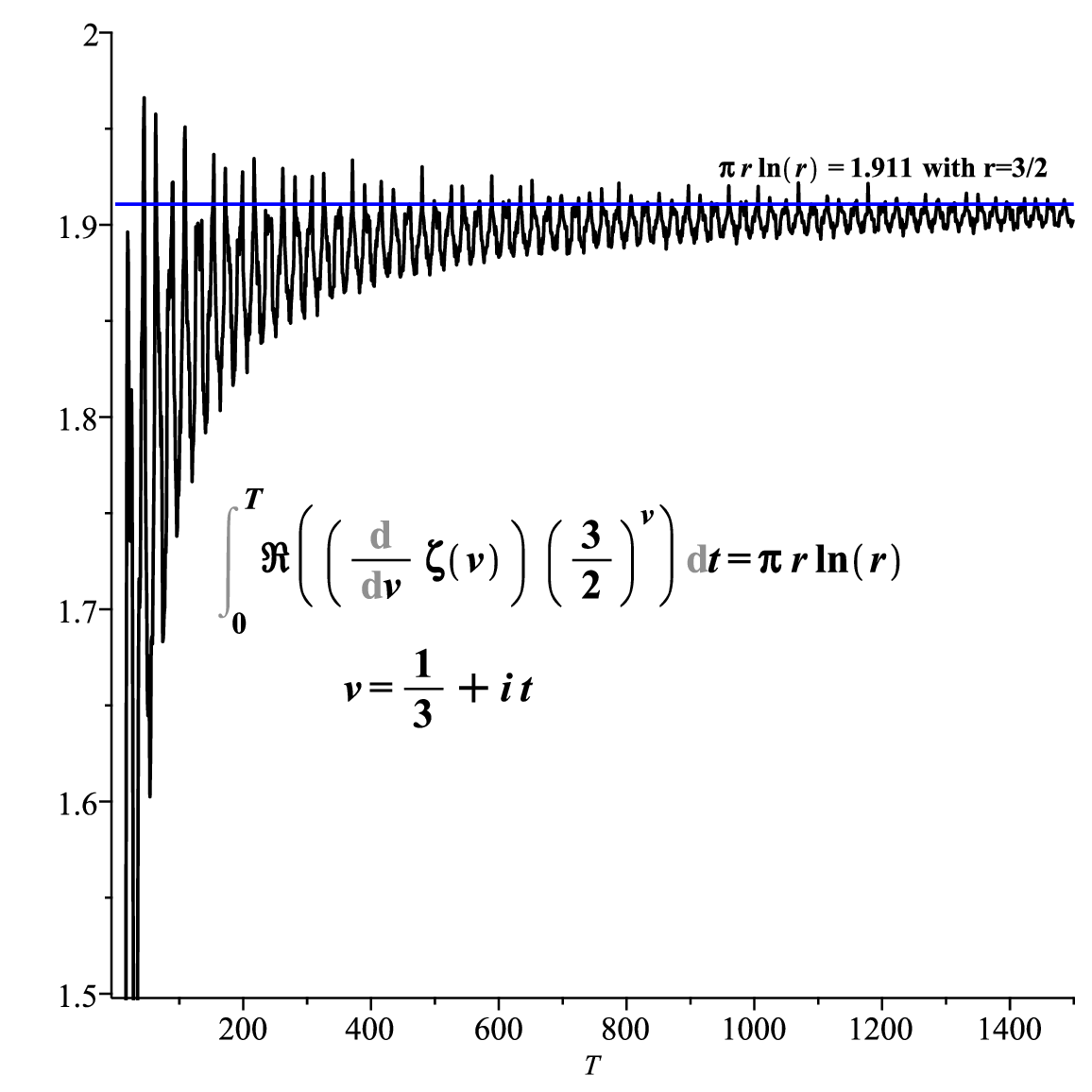}
\captionof{figure}{A test of \eqref{Rneqn2d} with $\sigma=1/3,~r=3/2$.}
\label{fig:Derivtest1} 
\end{minipage}%
\hfill
\begin{minipage} {0.40\textwidth} 
\centering
{\includegraphics[width=1.0\textwidth,height=1.0\textwidth]{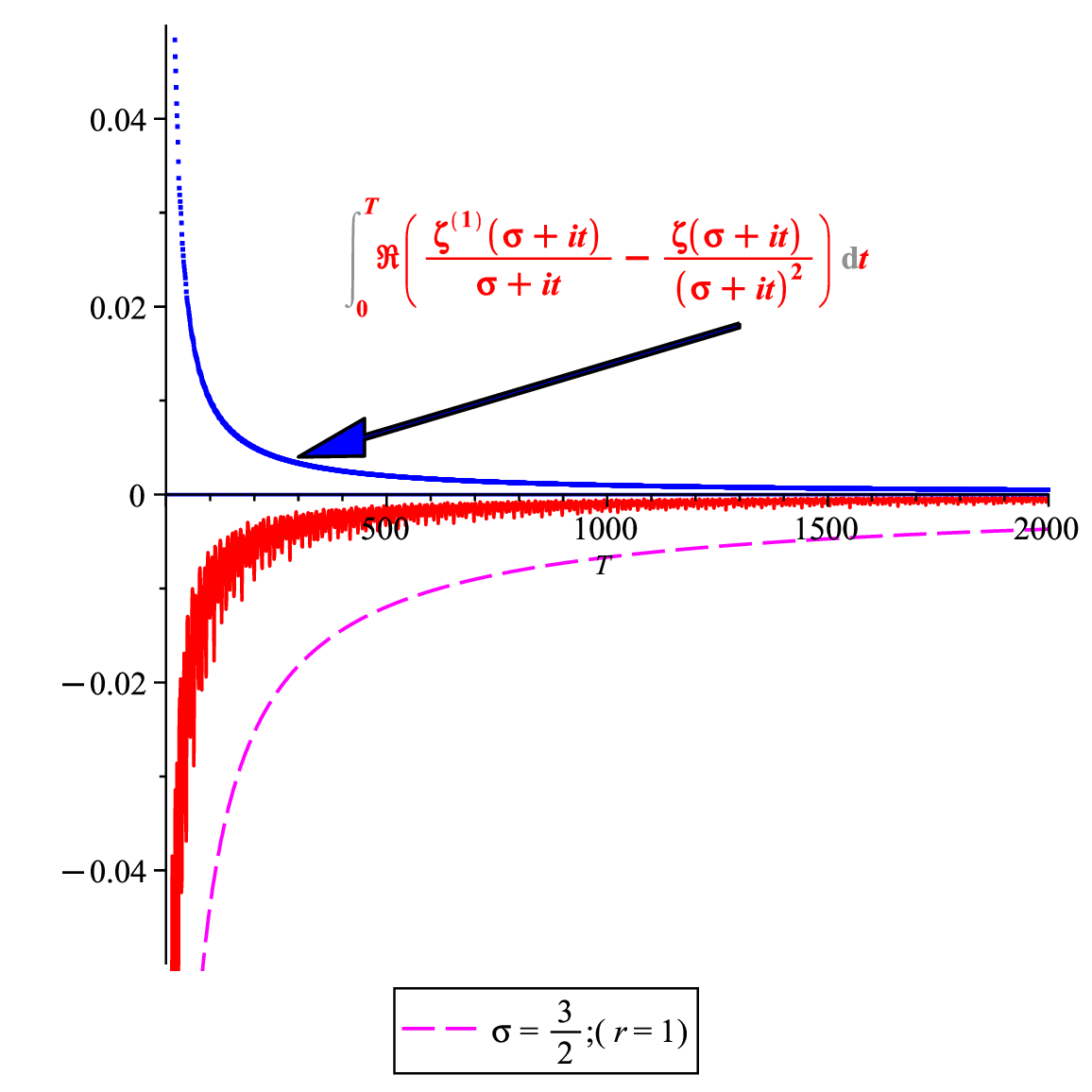}}
\caption{This Figure shows the partial sum of both the difference of the two integrals (red jagged curve) and the second integrand term itself (connected by the arrow), as well as the Cesàro estimate (dashed curve), all of which are tending to zero.}
\label{fig:DiffSig1p5r1p0}
\end{minipage}%
\end{figure}
\subsection{Derivatives}
\begin{figure} [H] 
\center
\begin{subfigure}   
[{Test of \eqref{Jd1} (i.e., $\sigma>1$ and $r\neq n$), calculated as a difference of two integrals.}]
{\includegraphics[width=0.40\textwidth,height=0.40\textwidth]{{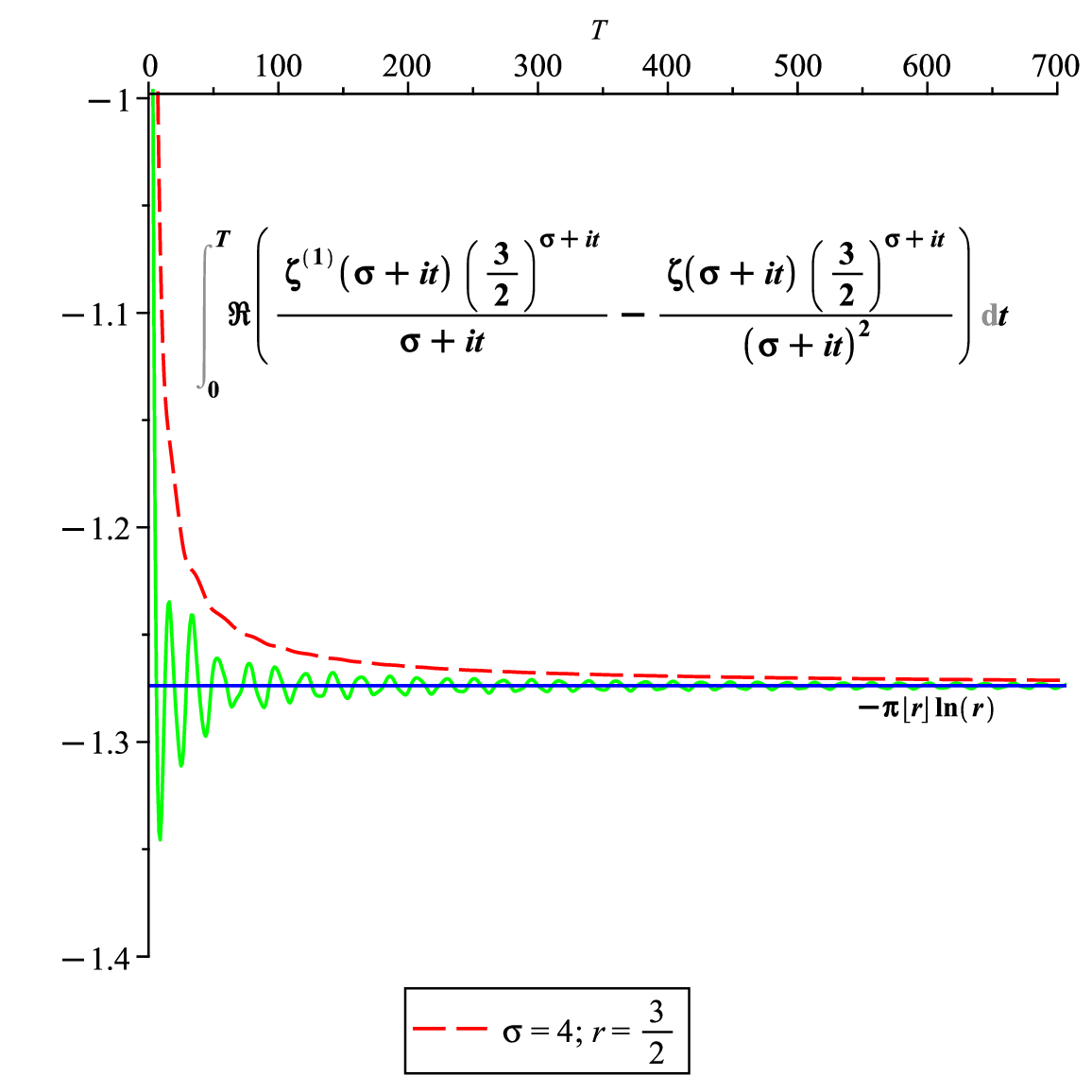}}
}
\label{fig:DiffSig3p0_r2p0}
\end{subfigure}
\hfill
\begin{subfigure}
[{Test of \eqref{Jdn} (i.e., $\sigma>1$ and $r=n$), calculated as a difference of two integrals.}]
{
\includegraphics[width=0.40\textwidth,height=.40\textwidth]{{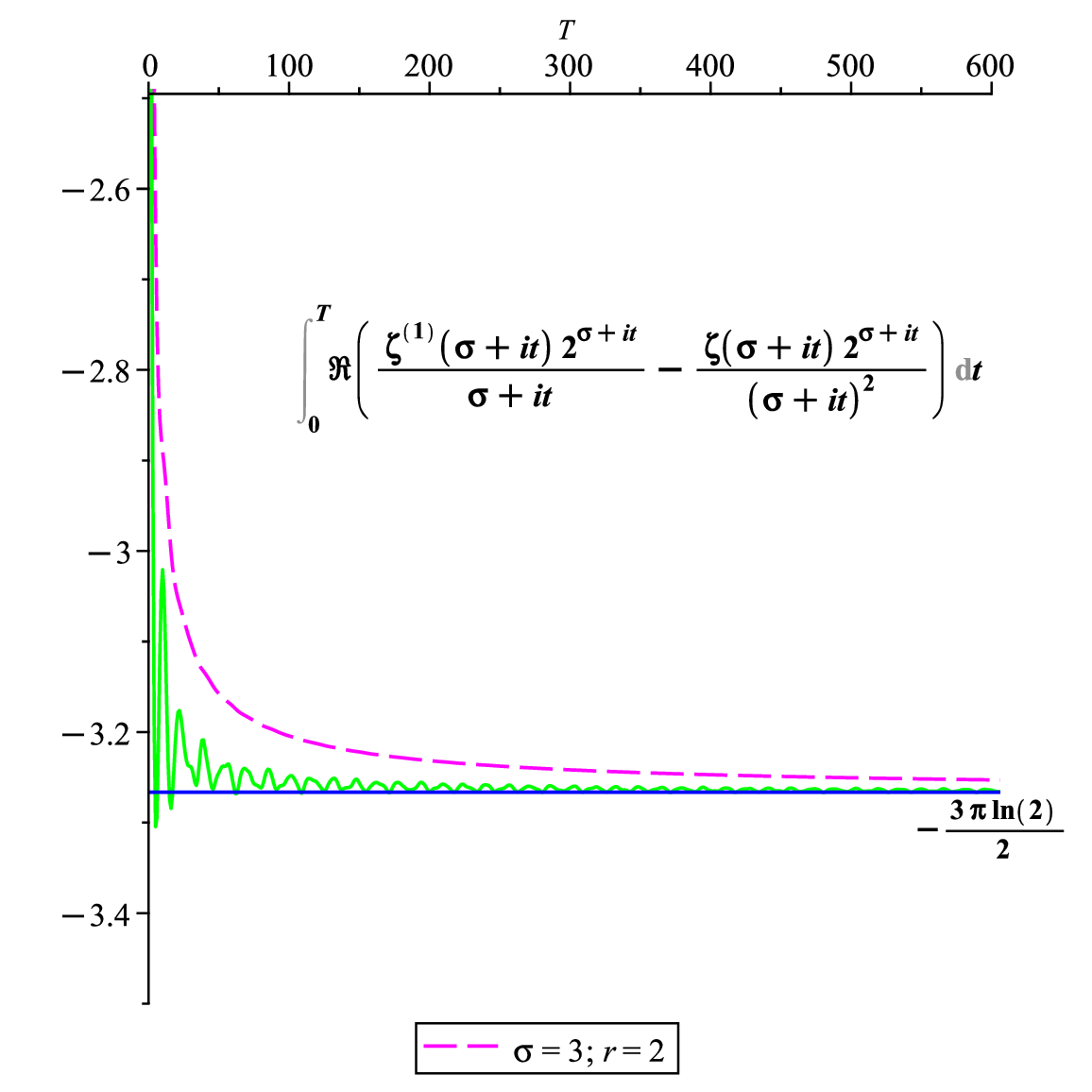}} 
}
\label{fig:DiffSig4p0_r2p0}
\end{subfigure}
\caption{This Figure tests \eqref{Jd1} and \eqref{Jdn}, both valid for $\sigma>1$.}
\label{fig:Drvtest}
\end{figure}

To test the identities involving derivatives, consider Figure \ref{fig:Derivtest1}, a test of \eqref{Rneqn2d} showing convergence to the expected limit. Mindful of the ambiguity discussed with respect to \eqref{Jdn1}, we next consider Figure \ref{fig:DiffSig1p5r1p0} showing that the difference of the two integrals (when $\sigma>1$) is tending to zero as expected.  Figure \ref{fig:Drvtest} shows that \eqref{Jd1} and \eqref{Jdn} also appear to be approaching their respective limits as predicted.  
%
%
%

\subsection{Periodicity} \label{sec:periodicity}
\begin{figure} [h]
\centering
\begin{subfigure}  
[{The small down arrow marks the point corresponding to the fifteenth harmonic of the assumed period $\rho$. The two other arrows indicate the minimum of the respective Cesàro approximation when $\sigma=1/2$ and $\sigma=1$. The partial sum corresponding to the numerical integration is also shown crossing the line ( $-2.1\pi$ ) that marks the asymptotic limit.}]
{\includegraphics[width=0.4\textwidth,height=0.5\textwidth]{{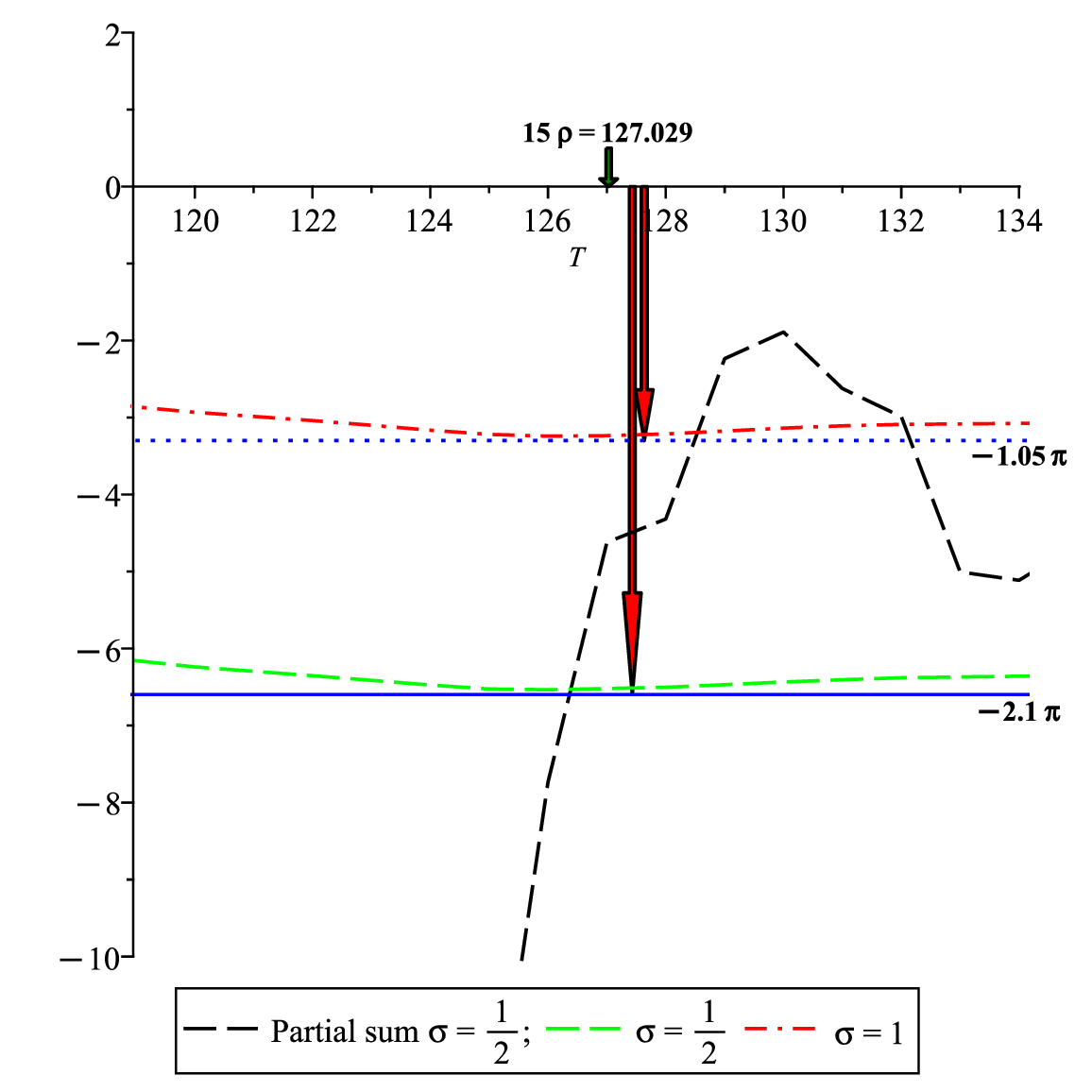}}
\label{fig:FirstCross2p1}
}
\end{subfigure}
\hfill
\begin{subfigure}
[{The small down arrows mark the points at a multiple of the harmonic in Figure \ref{fig:FirstCross2p1}. The two other arrows indicate the minimum of the respective Cesàro approximation when $\sigma=1/2$ and $\sigma=1$. The partial sum corresponding to the numerical integration up to the point $T$, is also shown crossing the lines ($-1.05\pi$) and ($-2.1\pi$ ) that mark the asymptotic limit for $\sigma=1$ and $\sigma=1/2$ respectively. When $\sigma=1/2$, this occurs at $T=1549.31$.}]
{
\includegraphics[width=0.50\textwidth,height=.5\textwidth]{{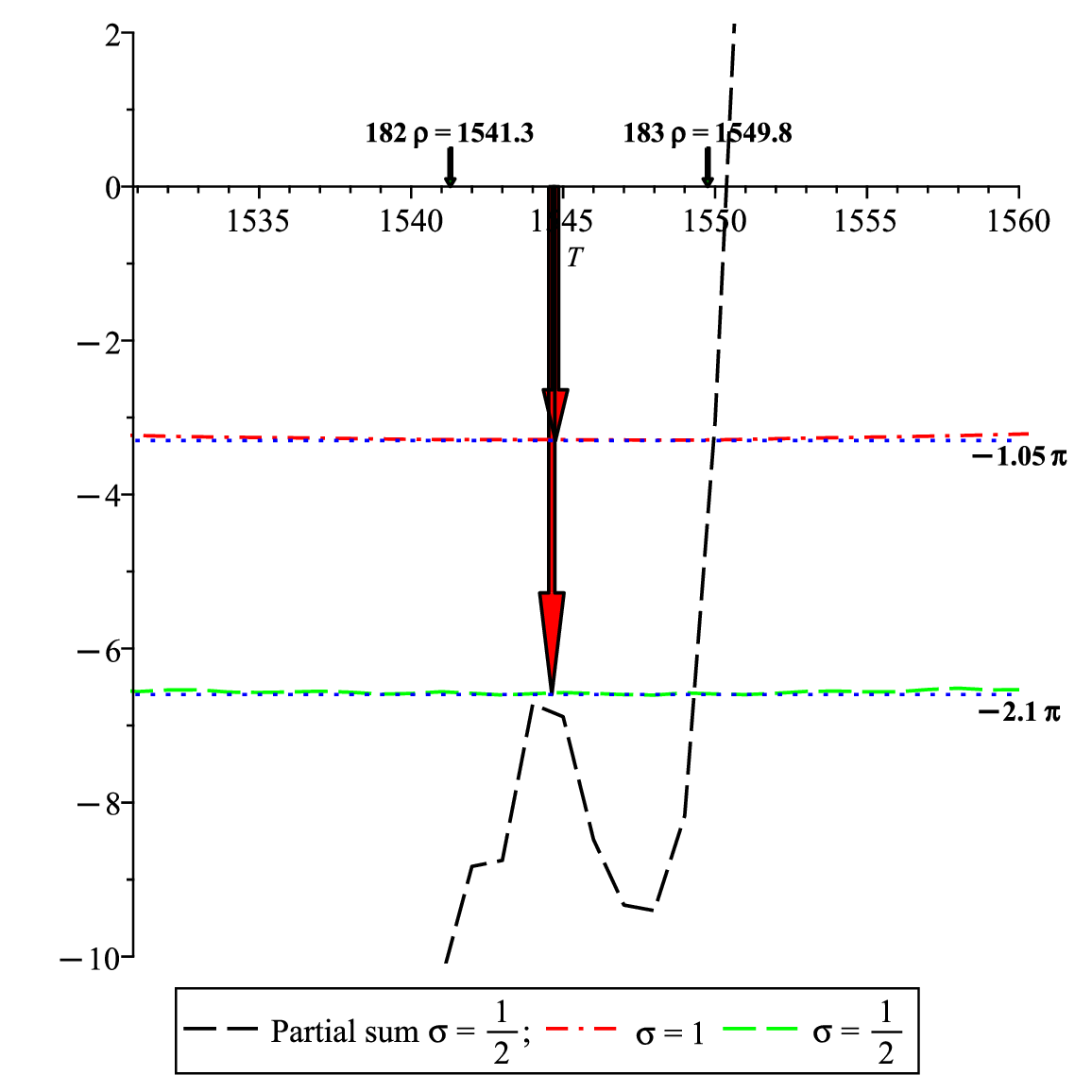}} \label{fig:Cross12}
}
\end{subfigure}
\caption{Focus on the Cesàro approximations to the function defined in \eqref{sighalf} with $r=2.1$ and various $0<\sigma\leq 1$. The tested cases $\sigma=1/3, ~ 3/4$ and $\sigma=1/2$ are indistinguishable within graphical resolution.}
\label{fig:Crossings}
\end{figure}

The periodicity of the estimates (e.g. Figure \ref{fig:3Sigs}), suggests the existence of a relationship with the periodicity of the integrand trigonometric functions, with a period being integral multiples of a quantity $\rho$ that satisfies the identity $2m\pi=\rho\ln(2.1)$, that is, if $m=1$, $\rho=8.47$ -- see \eqref{RsinCos}. The arrows in Figure \ref{fig:3Sigs} indicate an approximate observed periodicity of the Cesàro estimate equal to $15\rho$, with a small phase shift as $T$ increases. At a higher multiple of the initial period ($\sim180\rho$) the minimum appears to be shifted slightly. This is examined in more detail in Figure \ref{fig:Crossings}.
Because the minimum of the Cesàro estimate also appears to numerically coincide with the asymptotic ($T\rightarrow\infty$) value of the integral ($-2.1\pi=-6.597$), Figure \ref{fig:FirstCross2p1} isolates the region $T\approx 15\rho$ where it is shown that the near coincidence with the asymptotic result $-2.1\pi$ is close to, but not exactly in phase with the suspected periodicity within the accuracy of the numerical estimates used. Encouraged to speculate further by these near coincidences, consider Figure \ref{fig:Cross12} where the same estimates are presented when $T\approx 183\rho$. Here we find a very near-coincidence between the points where the partial sum crosses the asymptotic line ($-2.1\pi$) and the location of the ``harmonic" $183\rho$. For comparison, the location of the $182^{nd}$ harmonic is indicated by its own arrow and we note the the minimum of the Cesàro estimate is bounded by those two points. 

\begin{figure}
\center
%
%
\includegraphics[width=0.3\textwidth,height=0.3\textwidth]{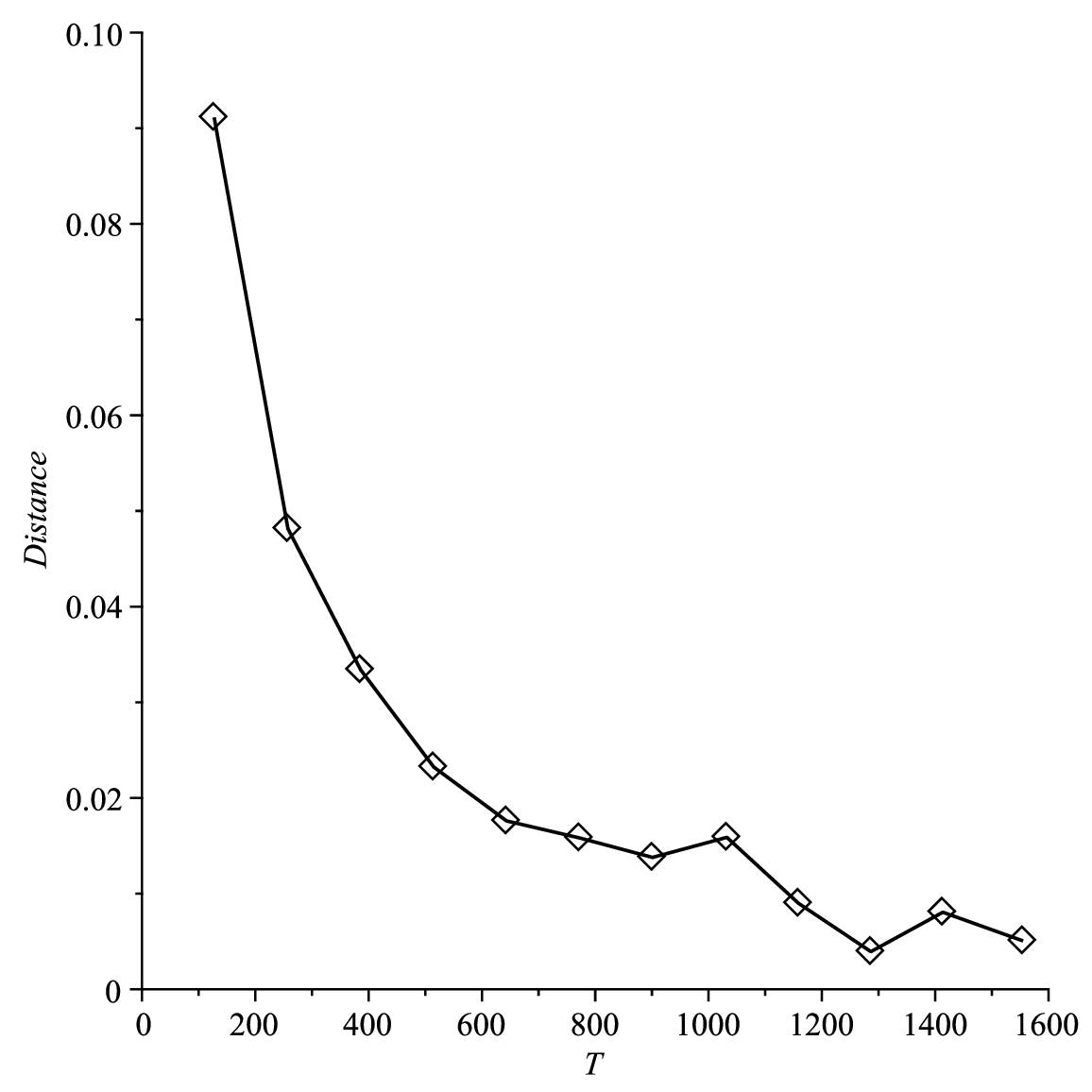} 
\captionof{figure}{For each minimum marked in Figure \ref{fig:3Sigs} this shows the distance between the Cesàro estimate and the asymptotic line $-2.1\pi$. See down arrows in Figure \ref{fig:Crossings}.}
\label{fig:XDists2p1}
\end{figure}

Although the Cesàro estimate appears to coincide with the asymptotic line in all these Figures, in fact careful study shows that it always lies a small distance above that line and by reasonable interpolation, it is possible to determine the point and distance of closest approach. These points are marked by large down arrows in Figure \ref{fig:Crossings} and presented quantitatively in Figure \ref{fig:XDists2p1} where we can see that the distance of closest approach decreases as $T$ increases, as expected of a converging approximation. Furthermore, by reasonable interpolation of the partial sums as they cross the asymptotic line, it is possible to test the possibility that the integral itself is periodic in integer multiples of the variable $\rho$; although the period between successive intersections was found to be within a few percent of the integral ratio $1:2:3\cdots$, the values were never close enough to perfect integers to validate such a hypothesis, and, after examining several variations with different values of $r$, it seems that any postulated universal multiplier increases with $r$ and does not appear to exist. 

\section{Correlations} \label{sec:Correlat}

As has been seen, the finite version of the integral $Z^{\prime}(\sigma,r)$ is oscillatory with a periodicity that is approximately constant if $r\neq n$, suggesting that segments of the integral itself are periodic and therefore correlated. Since the integral reflects the oscillatory nature of the integrand, it is worthwhile to conjecture that the integrand itself possesses corresponding segments with periodically correlated values. This could be a reflection of either (or both) of the cosine or absolute value factors appearing in \eqref{Rz}.

Specifically, if, as observed (see e.g., Figure \ref{fig:fig13b}), a partial sum and therefore its corresponding proper integral itself, is periodic, that is, in general, if
\begin{equation}
\int_{L_{1}}^{L_{2}}f \left(t \right)d t \approx 
\int_{L_{1}+\rho}^{L_{2}+\rho}f \left(t \right)d t,\hspace{50pt} \forall L_{1},L_{2},
\label{Gint}
\end{equation}
then
\begin{equation}
\int_{L_{1}}^{L_{2}}f \left(t +\rho \right)-f \left(t \right)d t
\approx 0\,,
\label{Gint2}
\end{equation}
and one possibility is that $f(t)$ and $f(t+\rho)$ are correlated to some extent. With reference to \eqref{Rz}, one expects a considerable degree of periodicity attached to the cosine function; however, more intriguing is the possibility that a correlation exists between elements of $\zeta(\sigma+it)$ lying on corresponding offset (i.e., shifted) segments of the imaginary line $\sigma+it$.
%
%
\begin{figure} [h]
\centering
\begin{subfigure}  
[{This shows overlain segments of the argument of \eqref{Rz}, each shifted relative to $t=0$ by the amounts $\rho$ indicated. In all cases, $r=2.1$.}]
{\includegraphics[width=0.4\textwidth,height=0.5\textwidth]{{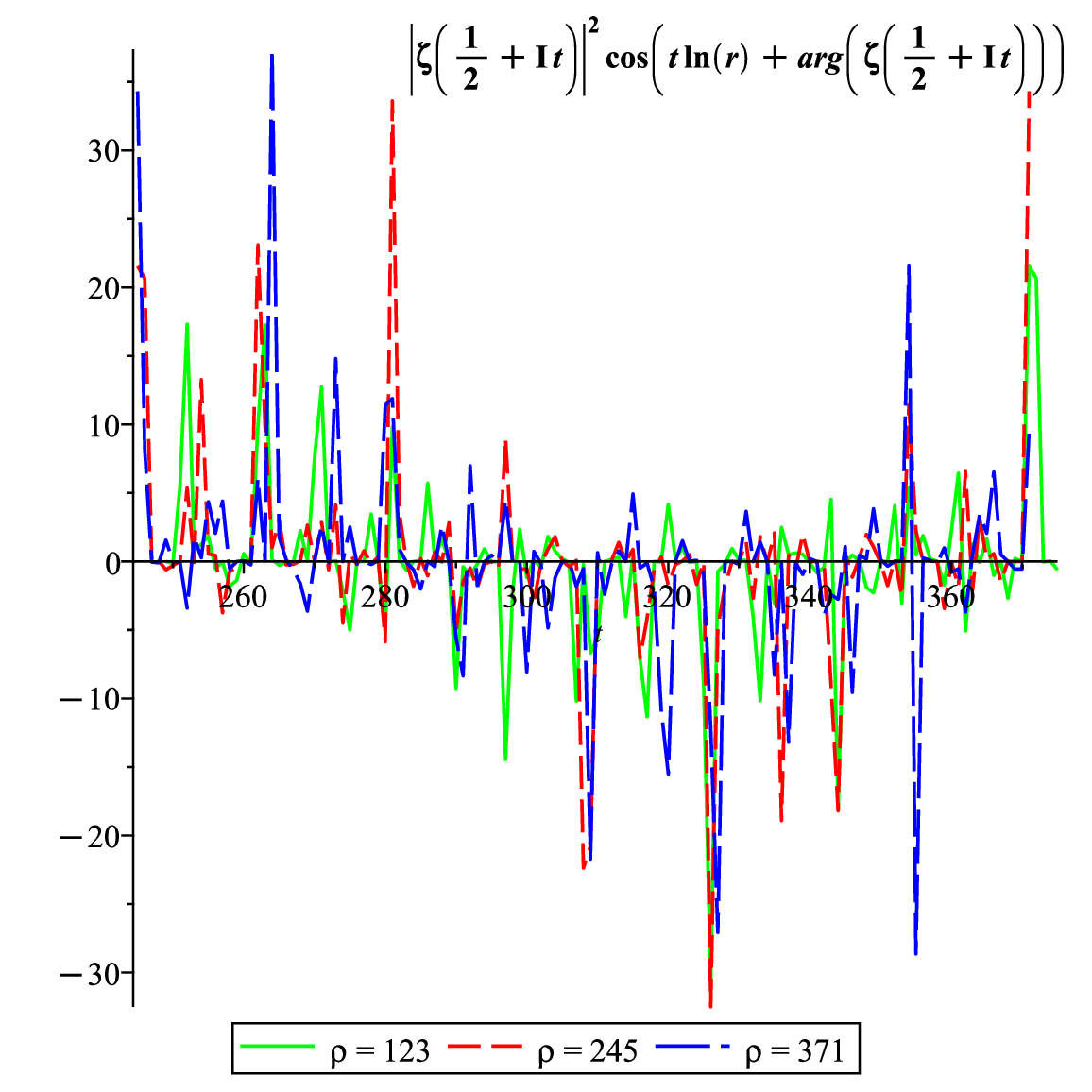}}
\label{fig:Shifted_Int}
}
\end{subfigure}
\hfill
\begin{subfigure}
[{This shows overlain segments of the function $\left|\zeta(1/2+it)\right|^2$, each shifted relative to $t=0$ by the amounts $\rho$ indicated.}]
{
\includegraphics[width=0.50\textwidth,height=.5\textwidth]{{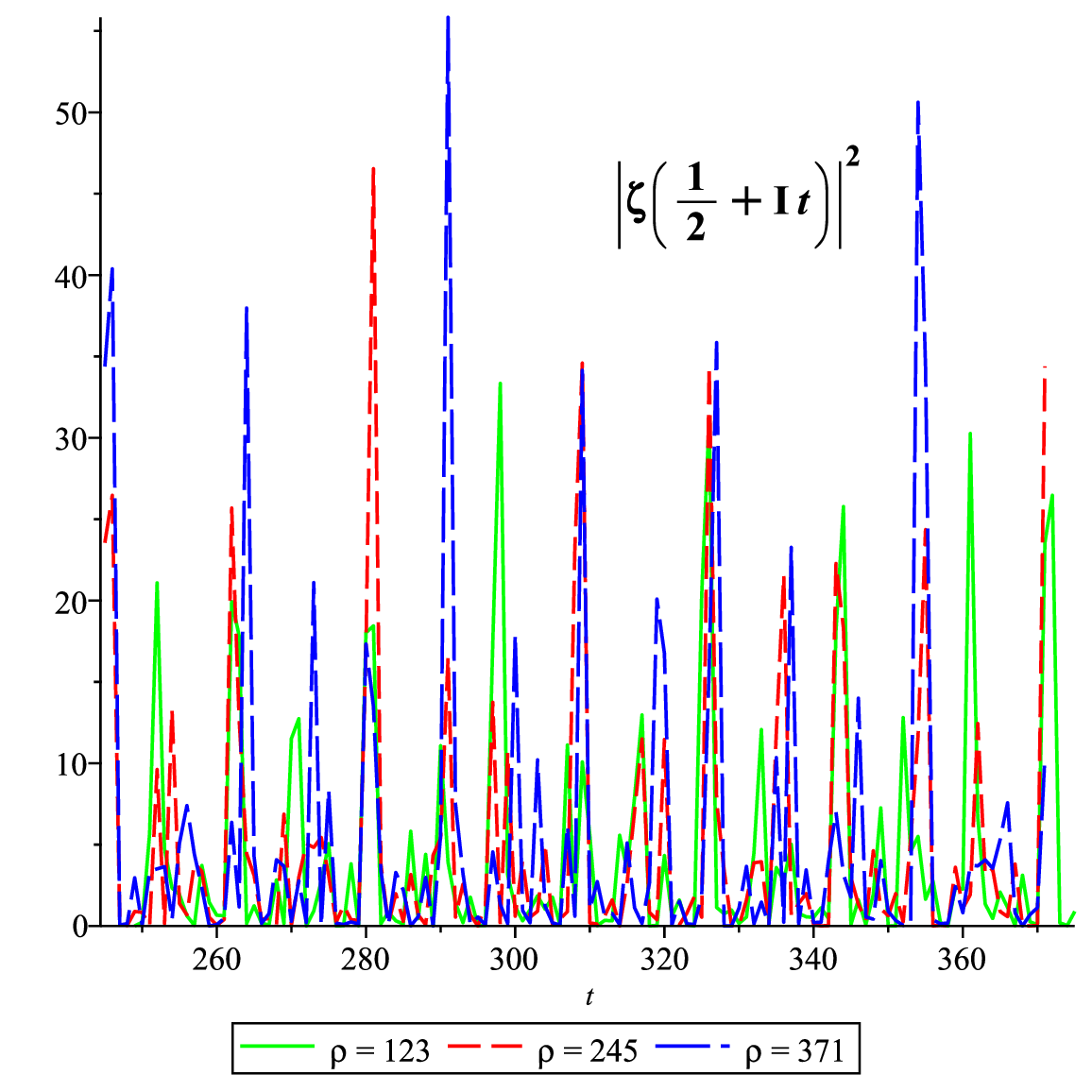}} \label{fig:ShiftAbs}
}
\end{subfigure}
\caption{This Figure overlays two quantities, each shifted by $t\rightarrow t+\rho$ shown.}
\label{fig:Shifts}
\end{figure} 
 
Again, inspired by Figures \ref{fig:3Sigs}, \ref{fig:MoreSigs2p1} and \ref{fig:fig13b}, Figure \ref{fig:Shifts}, which overlays three segments of the imaginary line, each of length 126 (see Figure \ref{fig:Crossings}) and each offset from $t=0$ by the quantity $\rho$ indicated, suggests that a correlation may, in fact, exist. This can be quantitatively explored by devising a coefficient to measure any correlation that may exist between segments of the $Z^{\prime}(\sigma,r)$ integrand as a function of the offset quantity $\rho$.

For any two continuous functions $f(t)$ and $g(t)$, we borrow from statistics, and define a correlation coefficient $\mathrm{Cor}(f,g)$ by

\begin{equation}
\mathrm{Cor}(f,g)=\frac{\mathrm{Cov}(f,g)}{\sqrt{\mathrm{Var}(f)\,\mathrm{Var}(g)}}
\end{equation}

where the ``Covariance" between two functions $f(t)$ and $g(t)$ is defined by
\begin{equation}
\mathrm{Cov}(f,g)\equiv E(fg)-E(f)E(g),
\end{equation}
the ``Variance" by
\begin{equation}
\mathrm{Var}(f)\equiv E(f^2)-E(f)^2
\end{equation}
in terms of the fundamental quantity ``Expectation" defined for any function $h(t)$ by
\begin{equation}
E(h)\equiv\frac{1}{L_{2}-L_{1}}\int_{L_{1}}^{L_{2}}h(t)\,dt\,.
\label{Expect}
\end{equation}
Effectively, the correlation coefficient measures the normalized (scaled between $\pm1$) correlation between the values of $f(t)$ and $g(t)$ averaged over the interval $\left[L_{1},L_{2}\right]$, with the statistical average operator replaced here by a numerical integration operator. If $\mathrm{Cov}(f,g)$ is {\it close} to $\pm1$, the functions $f(t)$ and $g(t)$ are strongly correlated/anticorrelated; if $\mathrm{Cov}(f,g)$ is {\it close} to zero, they are uncorrelated; otherwise the strength of the (anti)correlation becomes a subject of statistical study and opinion. Since we are not dealing here with statistical correlations between individual Cesàro summation elements, it is worth noting an adage that encapsulates varied opinion on this subject: {\it correlation coefficients between 0.5 and 0.7 indicate variables that are moderately correlated.} See also \cite{Stack233606}.
\begin{figure}
\noindent\begin{minipage} {0.47\textwidth} 
\centering
{
\includegraphics[width=1\textwidth,height=0.9\textwidth]{{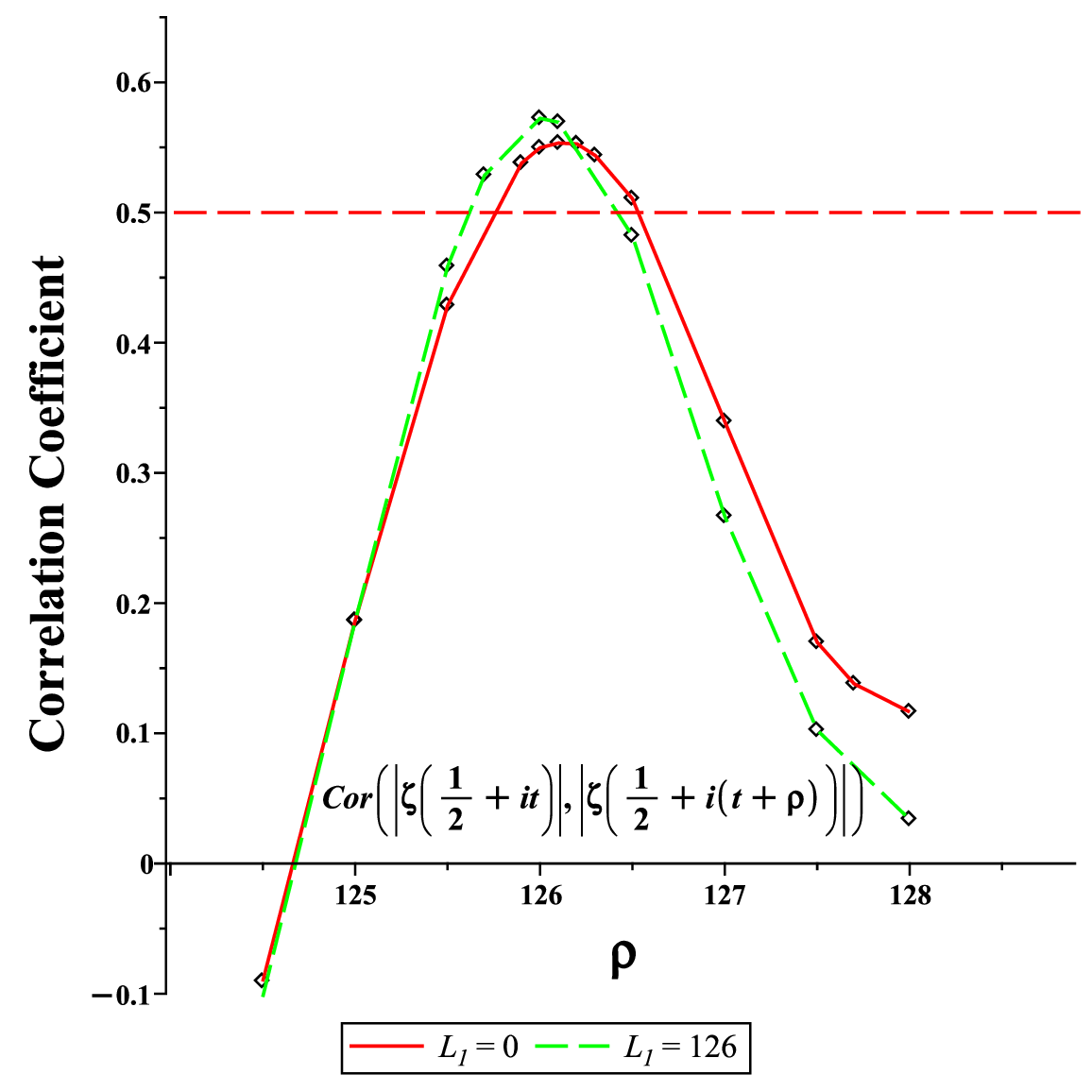}}
}
\caption{This Figure shows the correlation coefficient between $f(t)=\left|\zeta(1/2+it)\right|$ and $f(t)=\left|\zeta(1/2+i(t+\rho))\right|$ as a function of $\rho$, between offset intervals starting at $t=L_{1}$ as shown. Values above 0.5 are considered to be moderately correlated.}
\label{fig:AbsZetaCor}
\end{minipage}%
\hfill
\begin{minipage} {0.47\textwidth} 
\centering
{
\includegraphics[width=1\textwidth,height=0.8\textwidth]{{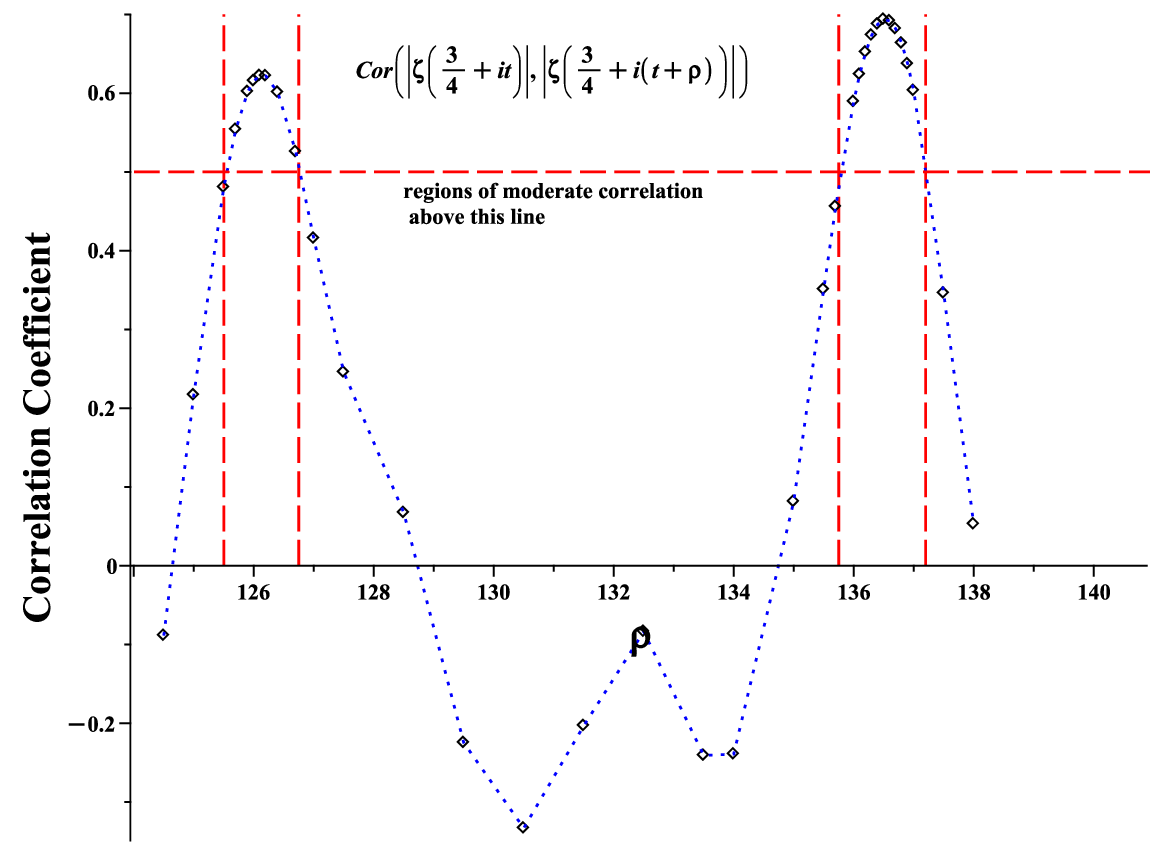}}
}
\caption{This Figure shows the correlation coefficient between $f(t)=\left|\zeta(3/4+it)\right|$ and $f(t)=\left|\zeta(3/4+i(t+\rho))\right|$ as a function of $\rho$.}
\label{fig:TwoShifts}
\end{minipage}%
\end{figure}

In the light of the observations discussed previously, and based on Figure \ref{fig:MoreSigs2p1}, in the following, set $L_{2}=L_{1}+126$. Figure \ref{fig:AbsZetaCor} shows the correlation coefficient $\mathrm{Cor}(f,g)$ between $f(t)=\left|\zeta(1/2+it)\right|$ and $g(t)=\left|\zeta(1/2+i(t+\rho))\right|$ as a function of $\rho$, evaluated over two contiguous pairs of segments, one starting at $L_{1}=0$, the other starting at $L_{1}=126$. This indicates that over a very narrow range of offset -- notably $\rho=126.1\pm 0.6$ -- there exist offset segments, between which the average value of $\zeta(1/2+it)$ compared between each segment, is moderately correlated. The fact that such contiguous segments exist and moderate correlation can be found for very precise values of the offset $\rho$, suggests that this is not a random statistical artifact. Therefore, it is reasonable to theorize that there exist other segments of the line $1/2+it$ among which the average value of the function $\left|\zeta(1/2+it)\right|$ is moderately auto-correlated \cite{Curran2023correlations}. This observation may in turn reflect on the distribution of zeros of the function $\zeta(1/2+it)$ (see  \cite{baluyot2023unconditional}). 

\begin{figure}
\includegraphics[width=1\textwidth,height=0.45\textwidth]{{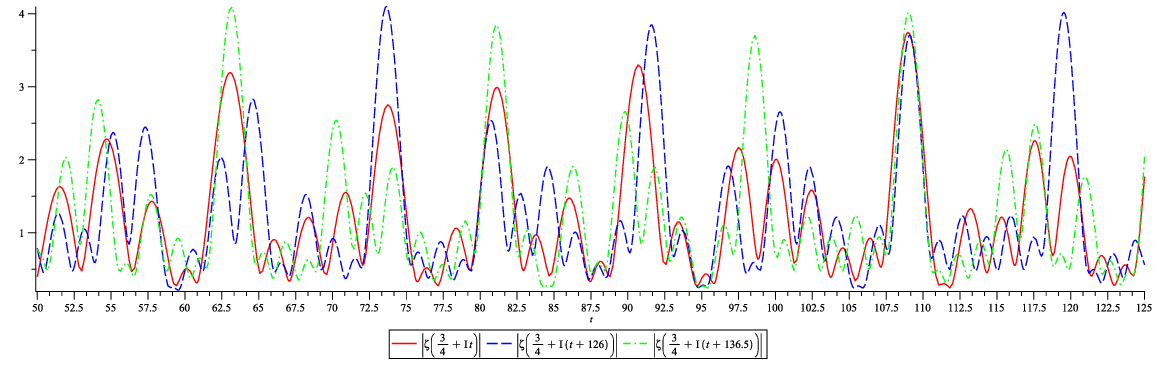}}
\caption{A comparison of $\left|\zeta(3/4+i(t+\rho))\right|$ with two different offsets over a small range of $t$. The corresponding comparisons for the real and imaginary parts of $\zeta(\sigma+it)$ are given in Appendix \ref{sec:TwoFigs}, Figures \ref{fig:ThreeshiftsRe} and \ref{fig:ThreeshiftsIm}}.   
\label{fig:Threeshifts}
\end{figure}

Since correlations between $\zeta(\sigma+it)$ and offset $\sigma$ have been reported elsewhere \cite[Section 6.4.1]{Milgram2019}, all of the above suggests that we evaluate the correlation coefficient for a wider range of $\sigma$ over a wider range of $\rho$. This is performed in Figure \ref{fig:TwoShifts} showing that moderate correlation exists over a much broader range of $\rho$ than the previous exercise would lead us to believe. In this Figure, we see that moderate correlation exists for offsets of $\rho=126.1\pm0.6$ and $\rho=136.5\pm0.7$ all of which suggests that there exists an underlying offset of $\rho\sim 10.1$ between correlated average values of $\zeta(\sigma+it)$. Qualitatively, the separation of peaks in Figure \ref{fig:Threeshifts} lends plausibility to this conjecture. The study of offset (i.e., shifted) moments of the zeta function involving integrals similar to those studied here also can be found in the literature (e.g. \cite{chandee2009correlation}, \cite{ng2022shifted}).

\section{Summary}

In this work, the functions $Z(\sigma,r)$ and $Z^{\prime}(\sigma,r)$ have been studied, both by analysis and numerical experimentation. As discovered both here and elsewhere, these functions possess the interesting property that for certain values of the parameter $r$, the value of the function depends continuously on its limiting value as $r$ approaches an integer $n$, but discontinuously at the point of discontinuity $r=n$. In particular, the traditional {\bf limiting} value of either function as $r\rightarrow n$ is different from its value {\bf at} $r=n$, demonstrating that both functions are discontinuous and the traditional definition $\underset{r\rightarrow n}\lim$ is invalid, but only at the limit point. With respect to the function $Z(\sigma,r)$, the value at the point of discontinuity $r=n$ lies naturally, rather than by decree, halfway between its value(s) as the limit point is approached from either direction. This property was verified numerically.

The derivative function $Z^{\prime}(\sigma,r)$ was also shown to have a similar property, notably that, for a  certain range of the parameter $\sigma$, it vanishes everywhere except at points $r=n$ where it is divergent. In this respect it represents one tine of the Dirac comb function, without the necessity of invoking limits or test functions. Again, this was verified numerically -- the numbers speak for themselves.

A very significant outcome of this study was the demonstration of the power of a modified Cesàro summation to approximate an improper integral and to verify analytic predictions. It is well-known (e.g., see Figures \ref{fig:Threeshifts} and \ref{fig:fig13a}) that the Riemann zeta function fluctuates wildly, but somehow the Cesàro summation tames the fluctuations and exposes properties that would otherwise have remained well-hidden. Certainly, without the Cesàro approximation, the quasi-periodic nature of $Z(\sigma,r)$ and the offset property of $\left|\zeta(1/2+it)\right|)$ would be unrecognized.

Having brought all these interesting properties to light, it is necessary to acknowledge that many questions remain, all of which are outside the scope of this study. Chief among them leads one to question ``Why is Cesàro approximation so effective?" Along the same lines, we observe that the Cesàro approximation appears to always approach its asymptotic limiting value from one direction -- is this significant and is there a reason? With respect to the fact that the modified Cesàro summation was originally introduced as a means of approximating an improper integral, it is necessary to recognize that although $T\sim4000$, the maximum used here, takes us a reasonable distance along the number line, it is still a long way to infinity. Thus there is no assurance of the validity of the underlying assumption -- the numerical observations presented will continue to behave as they have done here, as $T\rightarrow\infty$. 

When studying the periodic nature of the functions, it was first observed that the offset parameter $\rho=126$ yielded a reasonable suggestion of the existence of an underlying correlation, but that this choice is probably serendipitous and reflects a deeper periodicity of $\sim10$. If that is the case, the connection between a very clear Cesàro periodicity signal and a moderate underlying auto-correlation of the integrand function remains to be investigated. The underlying reason for the numerical value of any of the observed periodicities is not evident, and so it must be finally noted that the numerical experiments reported here merely scratch the surface of further study needed to verify these observations. 

\section{Caveats}
The  calculations and analysis in this work are based on two fundamental premises:
\begin{itemize}
\item{Glasser's Master Theorem (or, equivalently the translation of contour integrals) is valid only if the integrand functions vanish appropriately at the integration end-points;}
\item{The numerical integrations were performed initially using \cite[Maple]{Maple23} and later \cite[Mathematica]{Math23} when it was found that the latter executed orders of magnitude faster than the former and gave the same answers. However, the overall validity is dependent on the accuracy of algorithms built into these two computer codes on which our results are dependent.}
\end{itemize}

\section{Acknowledgement}
Professor Christophe Vignat (Tulane University) commented perceptively on the previous work \cite{milgram2024extension} and that comment led directly to the present study. He is also the author of the challenge Appendix \ref{sec:false}; we are grateful for his contribution.
\section{Funding}
All expenses associated with this work were funded by the authors.


\bibliographystyle{unsrt}
\bibliography{4p13_a.bbl}


\begin{appendix}

\section{Appendix: Cesàro Summation - a primer} \label{sec:Primer}

Cesàro summation \cite{WikiCesaro} assigns a value to an improper integral or sum that is not necessarily convergent in the usual sense. In the first case, the Cesàro summability of an integral involves the limit of the means of its partial integrals; in the second case it involves the limit of the means of its partial sums. Here, we study an improper integral in a different way, by first converting it into a sum, utilizing simple subdivision and numerical integration, and then employ Cesàro summation to evaluate the sum. In particular we employ the Cesàro  algorithm $(C,1)$ \cite{Emath}, defined as the limit, as $J$ tends to infinity, of the sequence of arithmetic means of the first $J$ partial sums of the series. Since we are here interested in improper integrals, some of which possess debatable convergence status, by transforming each integral into a corresponding sum and observing a large number of partial sums, we obtain, if convergent, a valid estimate of its value, and, if divergent, an accepted regularization, which may, or may not, tend to infinity.  To be specific, consider an integral
\begin{equation}
H(T)=\int_{0}^{T}h(t)dt
\label{Hint}
\end{equation}
where we are interested in the limit $T\rightarrow \infty$
. Subdivide the integral into small parts defined by
\begin{equation}
h \left(t_{j}\right) = \int_{t_{j}}^{t_{j}+\delta}h \left(t \right)d t
\label{htj}
\end{equation}
so that 
\begin{equation}
H \left(t_{m}\right) = 
\moverset{m}{\munderset{j =1}{\sum}}h \left(t_{j}\right)
\label{Htj}
\end{equation}
where $t_{1}=0$, $t_{J}+\delta=T$ and $H(T)=H(t_{J})$. As employed here, we set $\delta=1$ and evaluate each of the elements $h(t_{j})$ by numerical integration over the appropriate interval to arrive at an accurate numerical estimate of the value of the integral $H(t_{m})$ over some interval defined by the choice of $j\leq m$. We now form partial sums defined by
\begin{equation}
P_{k} = \moverset{k}{\munderset{m =1}{\sum}}H \left(t_{m}\right),\hspace{50pt} k\leq J,
\label{Psum}
\end{equation}
and form the running average of each partial sum
\begin{equation}
H _{n} = 
\frac{1}{n}\moverset{n}{\munderset{k =1}{\sum}}P_{k}\,.\hspace{50pt} n\leq J
\label{Hn}
\end{equation}

\begin{figure}[h] \label{fig:fig13}
\centering
\begin{subfigure} [{Individual integration elements $h(t_{j})$.}]
{
\includegraphics[width=0.45\textwidth,height=.45\textwidth]{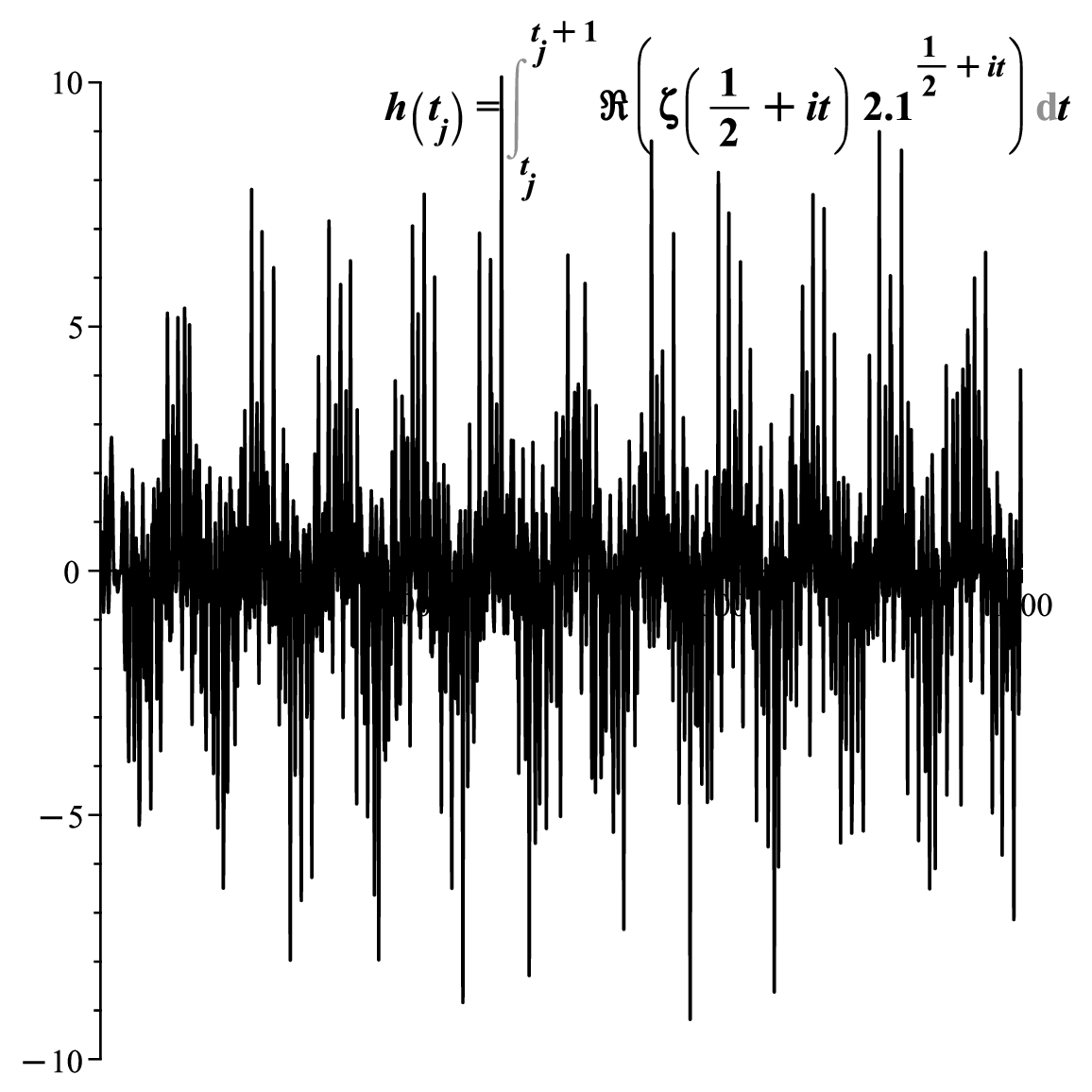} \label{fig:fig13a}
}
\end{subfigure}
\hfill
\begin{subfigure}[{The corresponding partial sums $P_{k}$. These also equate to the numerical value of the integral \eqref{RFt} as a function of its upper limit $T$. The asymptotic limit $-2.1\pi$ is also shown.}]
{
\includegraphics[width=0.45\textwidth,height=.45\textwidth]{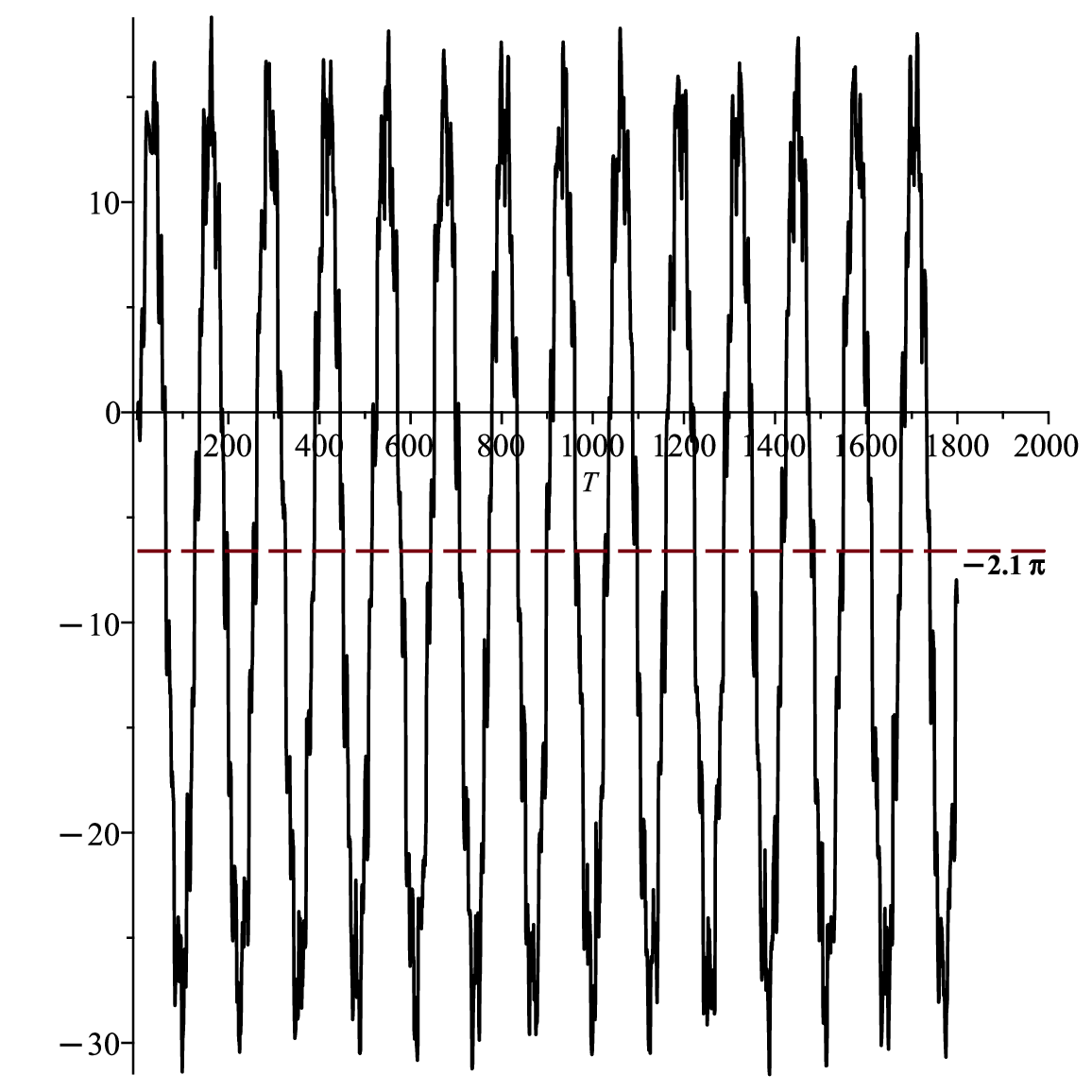} \label{fig:fig13b}
}
\end{subfigure}
\caption{A comparison of individual elements and partial sums. The average of the partial sums $P_{k}$ can be found in Figure \ref{fig:Z4(r)}}.
\end{figure}

If the sum is convergent, when $n= J$, $H_{J}$, the Cesàro sum, is an accurate numerical estimate \cite[VI.8]{VanDP} of the sum of the series \eqref{Htj}  and hence of the finite integral \eqref{Hint} and our interest focusses on $H_{J}$ as $J\rightarrow \infty$. The advantage of evaluating the integral in this manner is that it affords the analyst a simple means to study the properties of the integral as it converges to its upper limit, and, in the case of violently oscillating integrands, averaging the running partial sums tends to smooth any associated noise. If the integral diverges, it is still possible that the sum $H_{n}$ will converge to some finite value as $n\rightarrow J$  as $J$ itself increases, thereby yielding one possible regularization of a divergent series and its underlying divergent integral -- the sum of the series being defined (i.e., regularized) by the average of its partial sums. It is educational to consider an example examined in this study.

Figure \ref{fig:fig13a} demonstrates that the individual elements $h(t_{j})$ of the subdivision of the integral \eqref{sighalf} under consideration vary wildly as expected for the choice $\sigma=1/2$ and $r=2.1$. Figure \ref{fig:fig13b} shows that the partial sum of these elements varies considerably less violently than do the individual elements themselves. Since the partial sum equates to the numerical evaluation of the integral with varying upper limit, in this case we see that the finite integral periodically intercepts its own asymptotic value and is itself periodic. Finally, as presented previously in Figure \ref{fig:Z4(r)}, we see that the act of averaging the partial sums has completely washed out the noise, suggesting that the integral is dominated by the periodic components of its integrand. This is studied in Section \ref{sec:numtests}.
\begin{rem}
We note that it is the rare analyst who would pause to consider the average characteristics of the individual elements of a numerical integration and it is only through the use of Cesàro summation that one might consider such an examination.
\end{rem}

\section{Appendix: An enigmatic derivation} \label{sec:false}

As an alternative to the derivation presented in Section \ref{sec:method1}, here we consider an attempted independent  derivation of \eqref{Q1} and \eqref{Q2} by direct integration of an identity equivalent to the well-known \cite[Eq. (2.1.4)]{Titch2}, and convergent, integral representation 

\begin{align}
\frac{\zeta \! \left(\frac{3}{2}-i\,t \right)}{\frac{3}{2}-i\,t} &= 
\int_{1}^{\infty}\frac{{\lfloor x \rfloor}-x +\frac{1}{2}}{x^{\frac{5}{2}-i\,t}}d x +\frac{1}{\left(\frac{3}{2}-i\,t \right) \left(\frac{1}{2}-i\,t \right)}+\frac{1}{3-2\,i\,t}\\
&=\int_{1}^{\infty}\frac{{\lfloor x \rfloor}-x +\frac{1}{2}}{x^{\frac{5}{2}-i\,t}}d x +\frac{\frac{1}{2}+i\,t}{t^{2}+{1}/{4}}-\frac{\frac{3}{2}+i\,t}{2 \left(t^{2}+{9}/{4}\right)}\,.
\label{V1}
\end{align}
First, pre-multiply by $r^{-i\,t}$ and integrate \eqref{V1}, leading us to consider the integral

\begin{align} \nonumber
\int_{-\infty}^{\infty}\frac{\zeta \! \left(\frac{3}{2}-i\,t \right) r^{-i\,t}}{\frac{3}{2}-i\,t}d t
& = 
\int_{1}^{\infty}\frac{\left({\lfloor x \rfloor}-x +\frac{1}{2}\right) }{x^{\frac{5}{2}}}\int_{-\infty}^{\infty}\left(\frac{x}{r}\right)^{i\,t}d t\, d x\\&
 +\int_{-\infty}^{\infty}\frac{r^{-i\,t} \left(1/2+i\,t \right)}{t^{2}+1/4}d t -\frac{1}{2}\int_{-\infty}^{\infty}\frac{r^{-i\,t} \left(3/2+i\,t\right)}{t^{2}+9/4}d t 
\label{V4}
\end{align}
where the double integral operators in the first (right-hand side) term have been transposed. In a fashion similar to \eqref{J2s} and \eqref{G2}, with respect to the latter two integrals, we find

\begin{align} \nonumber
\int_{-\infty}^{\infty}\frac{r^{-i\,t} \left(\frac{1}{2}+i\,t \right)}{t^{2}+\frac{1}{4}}d t
& = 
\frac{1}{2}\int_{-\infty}^{\infty}\frac{\cos \left(t\,\ln \left(r \right)\right)+2\,t\,\sin \left(t\,\ln \left(r \right)\right) }{t^{2}+{1}/{4}}d t& \\ \nonumber
& = \frac{2\,\pi}{\sqrt{r}}  &\mathrm{if} ~r>1;\\ \nonumber
& =2\, \pi  &\mathrm{if} ~r=1;\\
& = 0  &\mathrm{if} ~0<r<1\,
\label{Jxa}
\end{align}
and
\begin{align} \nonumber
\frac{1 }{2}\int_{-\infty}^{\infty}\frac{r^{-i\,t} \left(\frac{3}{2}+i\,t \right)}{t^{2}+\frac{9}{4}}d t 
 &= 
\frac{ 1}{2}\int_{-\infty}^{\infty}\frac{\frac{3}{2}\,\cos \left(t\,\ln \left(r \right)\right)+t\,\sin \left(t\,\ln \left(r \right)\right) }{t^{2}+{9}/{4}}d t\\ \nonumber
 &= \frac{\pi}{r^{\frac{3}{2}}}  &\mathrm{if} ~r>1;\\ \nonumber
 & = \pi  &\mathrm{if} ~r=1;\\
 & = 0  &\mathrm{if} ~0<r<1\,.
\label{Jprint}
\end{align}
We now consider the first (right-hand side) term in \eqref{V4} and, making use of the identity (Maple, Mathematica) 

\begin{equation}
\int_{-\infty}^{\infty}\cos \! \left(t\,\ln \! \left(\frac{x}{r}\right)\right)d t
 = 2\,\pi \,\delta \! \left(\ln \! \left(\frac{x}{r}\right)\right)
\label{Cd}
\end{equation}

along with the change of variables $\ln(\frac{x}{r}):=x$, we obtain

\begin{tabular}{l l}
\\ 
\rule[-1ex]{0pt}{2.5ex} $ \int_{1}^{\infty}{x^{-\frac{5}{2}}}\left({\lfloor x \rfloor}-x +\frac{1}{2}\right)\, \int_{-\infty}^{\infty}\left(\frac{x}{r}\right)^{i\,t}d t d x $ &  \\ 

\rule[-1ex]{0pt}{5ex}  
$=\int_{1}^{\infty}{x^{-\frac{5}{2}}}\left({\lfloor x \rfloor}-x +\frac{1}{2}\right) \,\int_{-\infty}^{\infty}\cos \! \left(t\,\ln \! \left(\frac{x}{r}\right)\right)d t d x$ &  \\ 

\rule[-1ex]{0pt}{5ex} $=\,-{\pi}{r^{-\frac{3}{2}}}  \int_{-\ln \, \left({r}\right)}^{\infty}\left(2\,r\,{\mathrm e}^{x}-2\,{\lfloor r\,{\mathrm e}^{x}\rfloor}-1\right) {\mathrm e}^{-3x/2} \,\delta \! \left(x \right)d x$  &  \\ 

\rule[-1ex]{0pt}{5ex} $=\,-{\pi}{r^{-\frac{3}{2}}}  \left(2\,r -2\,{\lfloor r \rfloor}-1\right)$ & $\mathrm{if} ~r>1$ \\ 

\rule[-1ex]{0pt}{5ex} $= {\pi}{n^{-\frac{3}{2}}}  $ & $\mathrm{if} ~r=n,~n\geq 1;$ \\ 


\rule[-1ex]{0pt}{5ex}$=0$ & $ \mathrm{if} ~0<r<1\,.$ \\ 
\end{tabular}
\begin{equation} \label{Jcdefg}\end{equation}

\begin{rem} Nowhere in this derivation does $r=n$ suggest the existence of a special discontinuous case, and, in all cases, particularly \eqref{Jcdefg}, the case $r=n$ simply corresponds to the reduction $r\rightarrow n$ of the more general case corresponding to $r>1$. In contrast, see {\it Remark} \ref{sec:rgoeston}.
\end{rem}
Putting it all together, this derivation yields the prospective, but aberrant, identity
\begin{align} \nonumber
\int_{-\infty}^{\infty}\frac{\zeta \! \left(\frac{3}{2}-i\,t \right) }{\frac{3}{2}-i\,t}\,&r^{\frac{3}{2}-i\,t}d t\hspace{60pt} &
\\& \nonumber  = 2\,\pi \,{\lfloor r \rfloor} & \mathrm{if} ~r>1\,;
\\  \nonumber 
& =2\, \pi\,n & \mathrm{if} ~r=n,n\geq1;
\\  
& = 0 & \mathrm{if} ~0<r<1\,.
\label{Case123}
\end{align}

\begin{rem}
Note that 
\begin{itemize}
\item{}(i) the derivations \eqref{Q2} and \eqref{Case123} disagree when $r= n$;
\item{}(ii) the disagreement extends to the numerical evaluation presented in Figure \ref{fig:Z(r)}, and
\item{}
(iii) this derivation employs the traditional definition of the floor function (Remark \ref{rem:Mark1}).
\end{itemize}

\end{rem}
\section{Appendix: Two Figures} \label{sec:TwoFigs}
The following two figures present the real and imaginary components of the functions shown in Figure \ref{fig:Threeshifts}.
\begin{figure} [H]
\includegraphics[width=1\textwidth,height=0.45\textwidth]{{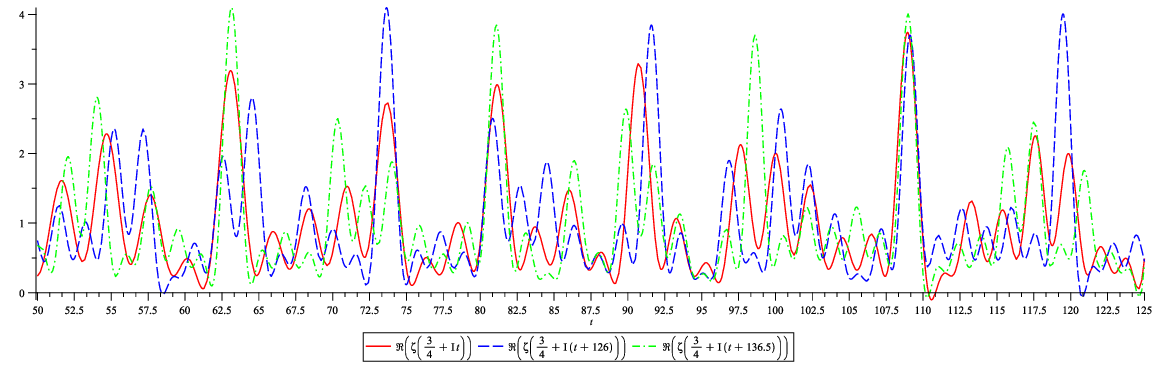}}
\caption{A comparison of $\Re\zeta(3/4+i(t+\rho))$ with two different offsets over a small range of $t$.}  
\label{fig:ThreeshiftsRe}
\end{figure}

\begin{figure} [H]
\includegraphics[width=1\textwidth,height=0.45\textwidth]{{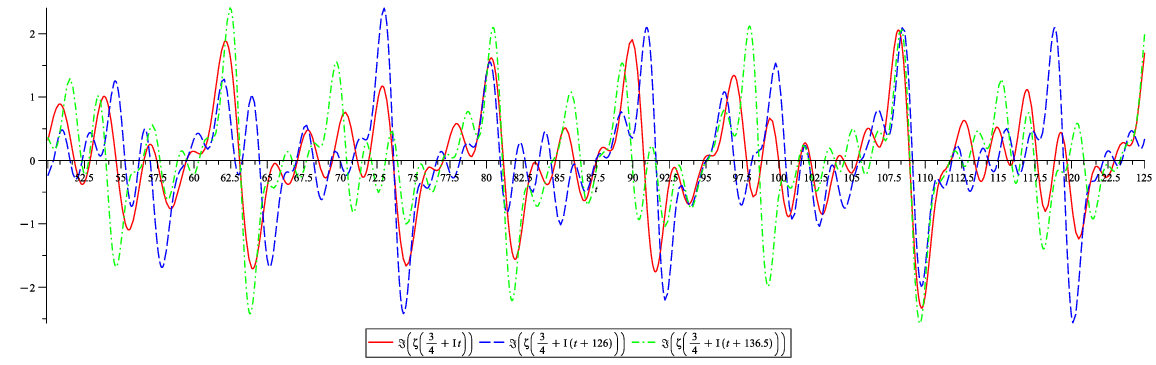}}
\caption{A comparison of $\Im\zeta(3/4+i(t+\rho))$ with two different offsets over a small range of $t$.}  
\label{fig:ThreeshiftsIm}
\end{figure}
\end{appendix}

\end{document}